\documentclass[reqno,11pt] {article}
\usepackage{amsmath, amsfonts, amssymb, amsthm,hyperref,floatrow}
\usepackage{graphicx, color,dsfont}
\usepackage{xcolor}
 \hypersetup{
    linktoc=page,
    linkcolor=red,          
    citecolor=blue,        
    filecolor=blue,      
    urlcolor=cyan,
    colorlinks=true           
}
\newcommand{\vct}[1]{{\mbox {\boldmath $#1$}}}

\newtheorem{theorem}{Theorem}

\newtheorem{lemma}[theorem]{Lemma}

\newtheorem{conjecture}[theorem]{Conjecture}

\def \P{ \mathbb P  }

\newenvironment{rem}[1][{\it Remark.}]{\begin{trivlist}
\item[\hskip \labelsep {\bfseries #1}]}{\end{trivlist}}

\definecolor{remi}{rgb}{1,0,0}
\usepackage{titlesec}
    \titleformat{\section}[hang]
        {\color{remi}{}\bfseries\filcenter\large}
        {\thesection.}
        {0.4em}
        {}[]

\DeclareMathSymbol{\leqslant}{\mathalpha}{AMSa}{"36} 
\DeclareMathSymbol{\geqslant}{\mathalpha}{AMSa}{"3E} 
\DeclareMathSymbol{\eset}{\mathalpha}{AMSb}{"3F}     
\renewcommand{\leq}{\;\leqslant\;}                   
\renewcommand{\geq}{\;\geqslant\;}                   

\newcommand{\R}{\mathbb{R}}

\title{Gaussian multiplicative chaos for symmetric isotropic matrices}

\author{{\sc Laurent Chevillard\thanks{e-mail: \texttt{laurent.chevillard@ens-lyon.fr}}} ,\, {\sc R\'{e}mi Rhodes\thanks{e-mail: \texttt{rhodes@ceremade.dauphine.fr}}}  ,\, {\sc Vincent Vargas\thanks{e-mail: \texttt{vargas@ceremade.dauphine.fr}}}
\\
{}\\
{\it $^*$Laboratoire de Physique de l'ENS Lyon, CNRS} \\
{\it 46, all\'ee d'Italie, F-69364 Lyon Cedex 07, France}  \\
\\
{}\\
{\it $^\dagger$ $^\ddagger$Universit\'e Paris-Dauphine, Ceremade, UMR 7564}\\
{\it Place de Mar\'{e}chal de Lattre de Tassigny}\\
{\it 75775 Paris  Cedex 16, France}\\
{}\\
}

\begin{document}
\maketitle
 \begin{abstract}
Motivated by isotropic fully developed turbulence, we define a theory of symmetric matrix valued isotropic Gaussian multiplicative chaos. Our construction extends the scalar theory developed by J.P. Kahane in 1985.
 \end{abstract}

\normalsize


\section{Introduction}
In the pioneering work \cite{Kah}, J.-P. Kahane introduced the theory of Gaussian multiplicative chaos. Given a metric space and a reference measure, Gaussian multiplicative chaos gives a mathematically rigorous definition to random measures defined as limits of measures with a lognormal density with respect to the reference measure. The main application of this theory was to define the Kolmogorov-Obhukov model of energy dissipation in a turbulent flow (see \cite{Kol,Man1}): in this context, the metric space is the euclidean space $\R^3$ equipped with the Lebesgue measure and the log density has logarithmic correlations. Since this seminal work, the theory of Gaussian multiplicative chaos has found many applications in a broad number of fields among which finance (\cite{BaKoMu,DuRoVa}) and 2-d quantum gravity (see \cite{Dav,KPZ} for the physics literature and \cite{BS,DS,RV,Rnew4} for the mathematics literature).

Three dimensional fluid turbulence is an archetypal out-of-equilibrium system in which energy is constantly injected at large scale and dissipated at the small viscous scales in a stationary manner. A statistical approach has been rapidly adopted in order to describe the complex multi-scale motions taking place in the flow. In the seminal work of Kolmogorov, known as the K41 theory \cite{Kol41,Fri95}, focusing on fully developed homogeneous and isotropic turbulence,  it is shown from the Navier-Stokes equations that energy is transferred from large to small scales at a constant rate, independently on viscosity: this is the fourth-fifth law. Further phenomenological extensions of this theory \cite{Kol,Fri95} took into account the peculiar statistical nature of the dissipation field that implies intermittent (or multifractal) corrections and probably more importantly, a long range correlated structure of velocity fluctuations. At this stage, scalar multiplicative chaos appears to be a good candidate to give a stochastic representation of the dissipation field, although nothing is said on energy transfers that ask for a stochastic model for the velocity field itself. Indeed, modern statistical studies underline the importance of defining a probabilistic model for the velocity field. Ideally, one looks for a field as close as possible to an invariant measure of the equations of motion (see for instance \cite{FoiMan01}). One of the first attempts in this direction is proposed in Ref. \cite{RobVar08} where the authors use a scalar multiplicative chaos to disturb an underlying Gaussian velocity field. A great success of this work is to propose an intermittent velocity field but unfortunately, the authors failed at proposing an incompressible dissipative velocity field, i.e. a field that respects the fourth-fifth law of Kolmogorov. One of the reasons is that the construction of the field does not include a basic mechanism of the Euler equations, namely the vorticity stretching phenomenon. This is the main novelty of the approach proposed in Ref. \cite{CheRob10}. One of the key steps of this construction is the introduction of the exponential of a Gaussian isotropic symmetrical random matrix in replacement of the scalar chaos used in Ref. \cite{RobVar08}. Heuristically, this symmetric matrix is reminiscent of the deformation field $\mathcal S$, i.e. the symmetric part of the velocity gradient tensor, that stretches vorticity $\vct{\omega} = \vct{\nabla} \wedge \textbf{u}$, where $\textbf{u}$ is the velocity field, according to the Euler equations:  $D\vct{\omega}/Dt =\mathcal S\vct{\omega}$, with $D/Dt = \partial/\partial t +\textbf{u}.\vct{\nabla}$ the Lagrangian derivative. It is easily seen that a Taylor development at short time $\tau$ of the former dynamics leads to a linear differential equation that can be solved using matrix exponentials of the initial deformation field $\tau\mathcal S(0)$. Then, logarithmic correlations and the free parameter $\gamma^2$ quantifying the level of intermittency are introduced by hands. A numerical investigation of the obtained velocity field shows indeed a mean energy transfer across scales.

As far as we know, there is no matrix valued theory of Gaussian multiplicative chaos that would be crucial in further understandings of the mechanisms at the origin of this energy transfer as observed numerically. The purpose of this work is thus to define such a theory for Gaussian symmetric and isotropic matrices. In the next section, we present the framework and the main results. Section 3 is devoted to the proofs of our main results. In the appendix, we gather general formulas which are useful in our proofs.



Notations: we denote by $M(\R^d)$ the set of measures on $\R^d$ and by $M_s(\R^d)$ the set of signed measures on $\R^d$. We denote by $\mathcal{S}( M_s(\R^d) )$ the set of symmetric matrices whose components belong to $M_s(\R^d)$. The $N$-dimensional identity matrix is denoted by $I_N$ and $P_N=(1)_{1\leq i,j\leq N}$ stands for the $N\times N$ matrix filled with the coefficient $1$ in each entry.

\section{Framework and main results}\label{eremat}
We first motivate the structure of our Gaussian matrix-valued random field. We remind that a random matrix $X$ is isotropic if for any real orthogonal matrix $O$, the matrices  $X$ and $OX ^tO$ have the same probability law (where ${}^tO$ denotes the transpose of the matrix $O$). 
If $N$ is an integer and if one considers a centered symmetric isotropic Gaussian $N\times N$-random matrix $(X_{i,j})_{1 \leq i,j \leq N}$, it takes on the following structure (see Lemma \ref{isotgauss}):
\begin{itemize}
\item the diagonal entries $(X_{1,1},\dots,X_{N,N})$ are independent of the off-diagonal entries $(X_{i,j})_{i<j}$,
\item the covariance matrix of the diagonal entries is given by $(1+c)\sigma^2I_N-c\sigma^2P_N$ where $c\in ]-1,\frac{1}{N-1}]$ and $\sigma^2\geq 0$,
\item the off-diagonal entries $(X_{i,j})_{i<j}$ are mutually independent with variance $\sigma^2\frac{1+c}{2}$.
\end{itemize}
Therefore, if one wishes to consider a general Gaussian field of symmetric isotropic matrices, the natural construction of such a field is to introduce a spatial structure preserving the above statistical structure. This is the main motivation for the construction of our field, which we describe now.

We introduce a probability space $(\Omega, \mathcal{F},P )$ and denote expectation by $E$. We want to define a homogeneous field of symmetric isotropic Gaussian matrices with logarithmic spatial correlations. The spatial correlation structure will be encoded by a kernel $K:\R^d\to\R$  of positive type of the form
$$K(x)=\gamma^2 \ln_+ \frac{L}{| x|}+g(x )$$ where $g$ is some continuous bounded function (in the sequel, we set $g(0)=m$) and $L>0$. Due to the divergence of this kernel at $x=0$, it is well known that the construction of such a field requires a cut-off approximation procedure to get rid of this singularity. Therefore, for $\epsilon>0$ (which stands in a way for the cut-off approximation rate), we introduce a covariance kernel $K_\epsilon:\R^d\to\R$ such that 
\begin{equation}\label{defKeps}
\sigma_{\epsilon}^2\stackrel{def}{=}K_\epsilon(0)=\gamma^2 (\ln \frac{L}{\epsilon}+1),\quad \text{and }\quad \sigma_{|y-x|}^2\stackrel{def}{=}K_\epsilon(x)=K(x) \text{ for all }|x|>\epsilon.
\end{equation}
Then we consider an integer $N\geq 2$ and $c \in ]-1,\frac{1}{N-1}]$. On this probability space,  for $\epsilon>0$, we consider a centered symmetric random matrix-valued Gaussian process $X^{\epsilon}(x)=(X^{\epsilon}_{i,j}(x))_{1 \leq i,j \leq N}$ indexed by $x\in\R^d$. We denote by
$$X^\epsilon_d(x)=(X^\epsilon_{1,1}(x),\dots,X^\epsilon_{N,N}(x))$$ the Gaussian vector made up of the diagonal entries of the matrix $X^\epsilon(x)$. We assume:
\begin{itemize}
\item the diagonal entries $(X^\epsilon_d(x))_{x\in\R^d}$ are independent of the off-diagonal entries $((X^\epsilon_{i,j}(x))_{i<j})_{x\in\R^d}$,
\item the covariance matrix kernel of the diagonal entries is given by 
$$E[\,{}^tX_d^\epsilon(x)X_d^\epsilon (y)]=\Big((1+c)I_N-cP_N\Big)K_\epsilon(x-y).$$ 
\item the off-diagonal entries $((X^\epsilon_{i,j}(x))_{i<j})_{x\in\R^d}$ are mutually independent, each of which with covariance kernel given by
$$ E[X^\epsilon_{i,j}(x)X^\epsilon_{i,j}(y)]=\frac{1+c}{2}K_\epsilon(x-y).$$
\end{itemize}
We also define
$$\bar{\sigma}_{\epsilon}^2\stackrel{def}{=}\frac{\sigma_{\epsilon}^2 (1+c)}{2}.$$

%


\begin{rem}\label{rem:kernel}
The canonical example of such a kernel $K$ is when it coincides with $\gamma^2 \ln  \frac{L}{| x|}$ for $x$ small enough. In dimension $1$ and $2$ we can even choose $K(x)= \gamma^2 \ln_+ \frac{L}{| x|}$. In dimension greater than $3$, we can use the constructions developed in \cite{Kah,RV1}: for examples of such kernels, see Appendix \ref{reminder}. Another approach is to use the convolution techniques developed in \cite{Rob}. This does not exactly fall into the framework set out above because the convoluted kernel depends on $\epsilon$ at all scales, i.e. for $|x-y|>\epsilon$. Nevertheless, this has  no significant influence on the forthcoming computations so that we also claim that our results remain valid for such regularization procedures.
\end{rem}

\begin{rem}
Note that the diagonal terms are independent if and only if $c=0$. In this case, the above structure coincides with the usual Gaussian Orthogonal Ensemble (GOE) \cite{Anderson,Mehta}. Note also that the boundary case $c=\frac{1}{N-1}$ corresponds to trace-free matrices.
\end{rem}

\begin{rem}
Application in turbulence. In the paper \cite{CheRob10}, the authors consider the following boundary case as a building block of their random velocity fields
\begin{equation*} 
\gamma^2= \frac{8}{3}\lambda^2, \; N=3, \; c=\frac{1}{N-1}=\frac{1}{2}
\end{equation*}
where $\lambda^2$ is found to fit experimental data for $\lambda^2 \approx 0.025$ \cite{CheRob10,CheCas12}. Here, the zero trace property is reminiscent of the incompressibility condition imposed on velocity fields.
\end{rem}

We want to study the convergence of the following random variable which lives in $\mathcal{S}( M_s(\R^d) )$
\begin{equation}  \label{eq:GaussChaos}
M^{\epsilon}(A)=\frac{1}{c_{\epsilon}}\int_A e^{X^{\epsilon}(x)} dx, \quad A \subset \R^d,
\end{equation}
where $c_{\epsilon}$ stands for a renormalization constant. From the scalar theory, we know that the constant $c_\epsilon$ is not trivial in order to avoid the blowing up of the above matrix integral. We will show that we can choose $c_\epsilon$ so as to have $E[ M^{\epsilon}(A)]= |A| Id$ where $|A|$ is the Lebesgue measure of $A$. Unlike the scalar theory, the explicit form of such a $c_\epsilon$ is not straightforward due to non-commutativity of the framework. We will prove that the normalization constant $c_{\epsilon}$ has the following explicit form
\begin{equation*}
 c_{\epsilon}=\frac{1}{N} \frac{\Gamma(1/2)}{\Gamma(N/2)}(1+c)^{(N-1)/2} \sigma_{\epsilon}^{N-1} e^{\frac{\sigma_{\epsilon}^2}{2}}
\end{equation*}
where $\Gamma$ stands for the usual Gamma function.

\begin{theorem}\label{Theo:Ncase}
Let $0<\gamma^2<d$. Then there exists a random matrix measure $M$ which lives in $\mathcal{S}( M_s(\R^d) )$ and such that for all bounded $A \subset \R^d$
\begin{equation*}
E[ {\rm tr} (M^{\epsilon}(A)-M(A))^2] \underset{\epsilon \to 0}{\rightarrow} 0.
\end{equation*}
We also have the following asymptotic structure
\begin{equation}
E[ {\rm tr} \left ( M(B(0,\ell)  )^2\right )]   
\underset{\ell \to 0}{\sim}  N^2 V_N \frac{\Gamma(N/2) e^{\gamma^2 \ln L +m}}{(1+c)^{(N-1)/2} \Gamma(1/2)} \frac{\ell^{2d-\gamma^2}} {(\gamma^2 \ln \frac{1}{\ell})^{(N-1)/2}}  \label{firstequivalent}
\end{equation}
with $V_N= \int_{|v|, |u| \leq 1} \frac{dudv}{|v-u|^{\gamma^2}}$. Furthermore, we get the following limit for every integer $k \geq 2$ such that $k<\frac{2d}{\gamma^2}$
\begin{equation} \label{secondequivalent}
\frac{\ln  E  [  {\rm tr} \left ( M( B(0,\ell)  )^k \right )   ] }{\ln \ell}  \underset{\ell \to 0}{\rightarrow}  \zeta(k)
\end{equation}
where $\zeta(k)=\big(d+\frac{\gamma^2}{2} \big)k-\frac{\gamma^2}{2} k^2 $.
\end{theorem}
 
Note that it would be interesting to prove that this matrix-valued Gaussian multiplicative chaos admits a phase transition as in the scalar case, which is likely to occur at $\gamma^2=2d$.

\begin{conjecture} 
Let $0<\gamma^2<d$. The power law spectrum of  $M$ is given by the following expression: for all $q \in ]0,\frac{2d}{\gamma^2}[$, $\forall \ell\in  (0,1]$,
$$E  [  {\rm tr} \left ( M( B(0,\ell)  )^q  \right )  ]\simeq C_q\ell^{\zeta(q)}(-\ln\ell)^{\frac{(q-1)(1-N)}{2}},$$ where $C_q>0$ is a constant and the structure exponent is given by
$$ \zeta(q)=\big(d+\frac{\gamma^2}{2} \big)q-\frac{\gamma^2}{2} q^2.$$
We give in the appendix a heuristic derivation of the above equivalent. If this conjecture is true, this would show that  noncommutativity yields an extra log factor in the power-law spectrum of $M$.
\end{conjecture}

\begin{rem}
Note that one can define a notion of ``metric" (actually a measure) through the quantity $$A\in \mathcal{B}(\R^d)\mapsto {\rm tr} M( A).$$
Therefore we can define the  notion of Hausdorff dimension associated to this ``metric" (see \cite{BS,DS,RV}). It would be interesting to prove a corresponding KPZ formula and relate it with a KPZ framework . 
\end{rem}

\section{Proofs of the $N$-dimensional case}
Let us first mention that several results about isotropic matrices and related computations are gathered in the appendix and will be used throughout this section.

\subsection{Joint law of the eigenvalues of Gaussian isotropic matrices}
We consider a symmetric random matrix $X=(X_{i,j})_{1 \leq i,j \leq N}$ made up  of centered Gaussian variables with the following covariance structure: the off-diagonal terms $(X_{i,j})_{i<j}$ are i.i.d. with variance $\sigma^2$. 
The diagonal term $(X_{1,1}, \cdots , X_{N,N})$ is independent from the off-diagonal and it has the following covariance structure
\begin{equation*}
K_N=( E[X_{i,i}X_{j,j}] )_{1 \leq i,j \leq N}=(1+c)\sigma_d^2 I_{N}-c\sigma_d^2 P_N
\end{equation*} 
where $I_N$ is the identity matrix, $P_N=(1)_{i,j}$ and $c \in ]-1,\frac{1}{N-1}[$. By noting that $P_N^2=N P_N$, we get the following inverse for $K$ if $c \not = \frac{1}{N-1}$
\begin{equation*}
K_N^{-1}=\frac{1}{\sigma_d^2(1+c)}I_{N}+\frac{c}{\sigma_d^2(1+c)}\frac{1}{(1+c(1-N))}P_N
\end{equation*}
The density of the random matrix, with respect to the Lebesgue measure $(dx_{i,j})_{i \leq j}$, is therefore given by
\begin{equation*}
 f\big( (x_{i,j})_{i \leq j} \big) = \frac{1}{Z_N} e^{  -\frac{1}{2\sigma_d^2(1+c)}\sum_{i=1}^{N}x_{i,i}^2-\frac{c}{2 \sigma_d^2(1+c)}\frac{1}{(1+c(1-N))}(\sum_{i=1}^{N}x_{i,i})^2 -\frac{1}{2\sigma^2}\sum_{i<j}x_{i,j}^2   }
 \end{equation*}
 where $$Z_N=(2 \pi)^{N(N+1)/4} \sigma_d^{N} \sigma^{N(N-1)/2} (1+c)^{(N-1)/2}\sqrt{1-(N-1)c}$$ is a normalization constant.
 
 Therefore if we have the following condition
 \begin{equation}\label{condisotrop}
 \sigma_d^{2}(1+c)=2 \sigma^2,
 \end{equation}
 as we have required in section \ref{eremat}, we can rewrite the above density in the following matrix form
 \begin{equation}\label{eq:1pointDensity}
 f( X )= \frac{1}{Z_N}e^{ -\frac{c}{2 \sigma_d^2(1+c)}\frac{1}{(1+c(1-N))}({\rm tr} X)^2 -\frac{1}{2\sigma_d^2(1+c)} {\rm tr} X^2   }
 \end{equation}
 with $Z_N=2^{N/2}\pi^{N(N+1)/4}\sigma_d^{N(N+1)/2}(1+c)^{(N-1)(N+2)/4}\sqrt{1+c(1-N)}$.
 This shows that the matrix is isotropic, namely that for any real orthogonal matrix $O$, the matrices  $X$ and $OX ^tO$ have the same probability law. Therefore by applying  \cite[Proposition 4.1.1, page 188]{Anderson}, we get the density of the unordered eigenvalues
 \begin{equation}\label{eq:densityvp}
 f( ( \lambda_{i} )_{1 \leq i \leq N} )=\frac{1}{\bar{Z}_N}e^{-\alpha(\sum_{i=1}^{N} \lambda_{i})^2-\frac{1}{2\sigma_d^2(1+c)}\sum_{i=1}^{N} \lambda_{i}^2} \Pi_{i<j} |\lambda_j-\lambda_i|,
 \end{equation}
 where $\alpha=\frac{c}{2 \sigma_d^2(1+c)}\frac{1}{(1+c(1-N))}$ and $\bar{Z}_N=2^{N(N-1)/4} \frac{\rho(U_1(\R))^N N! }{\rho(U_N(\R))} Z_N$ (notations of  \cite{Anderson}). We remind that $\rho(U_N(\R))=2^{N/2} (2\pi)^{N(N+1)/4} \prod_{k=1}^N \frac{1}{\Gamma(k/2)}$ (see \cite[page 198]{Anderson}) and thus
 \begin{equation}\label{eq:renormprobaeigen}
 \bar{Z}_N=N! (2\pi)^{N/2} (\prod_{k=1}^N \frac{\Gamma(k/2)}{\Gamma(1/2)}) \sigma_d^{N(N+1)/2} (1+c)^{(N-1)(N+2)/4}\sqrt{1+c(1-N)}. 
 \end{equation}
 The isotropic condition (Eq. \ref{condisotrop}) ensures also that the collection of eigenvectors $(v_i)_{1\leq i\leq N}$ is independent of the eigenvalues $( \lambda_{i} )_{1 \leq i \leq N}$, and they are distributed uniformly on the unit sphere according to the Haar measure \cite[Corollary 2.5.4, page 53]{Anderson}.

\subsection{Computations of the renormalization}

We consider here isotropic symmetric matrices $X^{\epsilon}(x)=(X^{\epsilon}_{i,j}(x))_{1 \leq i,j \leq N}$ as defined in section \ref{eremat} and compute the renormalization of order 1, i.e.  the constant $c_{\epsilon}$ such that
\begin{equation*}  
E[e^{X^\epsilon(x)}]=c_{\epsilon} I_N =  \frac{E[  \text{tr } e^{X^\epsilon(x)}  ]} {N} I_N.
\end{equation*}
The isotropic nature of the matrices ensures the proportionality of the former expectation to the identity matrix $I_N$ . We want more precisely an equivalent of $c_{\epsilon}$ as $\epsilon \to 0$. We have
\begin{equation*}
c_{\epsilon} =  \frac{1}{\bar{Z}_N} \int_{\R^N} e^{\lambda_1}e^{-\alpha_{\epsilon}(\sum_{i=1}^{N} \lambda_{i})^2-\frac{1}{2 \sigma_{\epsilon}^2(1+c)}\sum_{i=1}^{N} \lambda_{i}^2} \Pi_{i<j} |\lambda_j-\lambda_i| d \lambda_1 \cdots d \lambda_N,
 \end{equation*}
 where $\alpha_{\epsilon}=\frac{c}{2 \sigma_{\epsilon}^2(1+c)}\frac{1}{(1+c(1-N))}$ and the normalization constant $\bar{Z}_N$ is given by Eq. \ref{eq:renormprobaeigen} with $\sigma_d^2=\sigma_{\epsilon}^2= \gamma^2 (\ln \frac{L}{\epsilon}+1)$.
 
  We set  $u_i=\frac{\lambda_i}{ \sigma_{\epsilon}}$ and therefore we get
 \begin{equation*}
c_\epsilon=  \frac{\sigma_{\epsilon}^{N(N+1)/2}}{\bar{Z}_N} \int_{\R^N} e^{ \sigma_{\epsilon} u_1}e^{-\alpha (\sum_{i=1}^{N} u_{i})^2-\frac{1}{2(1+c)}\sum_{i=1}^{N} u_{i}^2} \Pi_{i<j} |u_j-u_i| d u_1 \cdots d u_N,
 \end{equation*}
 where $\alpha=\frac{c}{2 (1+c)}\frac{1}{(1+c(1-N))}$. We thus introduce 
 \begin{equation*}
 \varphi(u_1, \cdots, u_N)=  \sigma_{\epsilon} u_1-\alpha (\sum_{i=1}^{N} u_{i})^2-\frac{1}{2(1+c)}\sum_{i=1}^{N} u_{i}^2
 \end{equation*}
 The function $\varphi$ is maximal for $u_1=S_{\epsilon} (1+2 \alpha(1+c)(N-1) ), \;   i \geq 2:  \: u_i= -2 \alpha S_{\epsilon} (1+c)$ with $S_{\epsilon}=\frac{\sigma_{\epsilon}}{\frac{1}{1+c}+ 2 \alpha N}$. We thus set 
 $u_1=v_1+S_{\epsilon} (1+2 \alpha(1+c)(N-1) ), \;   i \geq 2:  \: u_i=v_i -2 \alpha S_{\epsilon} (1+c)$ to get
\begin{align*}
c_\epsilon & =  \frac{\sigma_{\epsilon}^{N(N+1)/2}e^{\frac{\sigma_{\epsilon}^2}{2}}}{\bar{Z}_N} \int_{\R^N}  e^{-\alpha (\sum_{i=1}^{N} v_{i})^2-\frac{1}{2(1+c)}\sum_{i=1}^{N} v_{i}^2}  \Pi_{2 \leq i}| v_1-v_i+(1+c)\sigma_\epsilon  |   \\
& \times  \Pi_{2 \leq i<j} |v_j-v_i| d v_1 \cdots d v_N.
 \end{align*} 
Therefore, we get the following equivalent by using the Laplace method
\begin{equation*}
c_\epsilon  \underset{\epsilon \to 0}{\sim}  \frac{\sigma_{\epsilon}^{N(N+1)/2}   (1+c)^{N-1} \sigma_{\epsilon}^{N-1} e^{\frac{\sigma_{\epsilon}^2}{2}}}  {\bar{Z}_N} \int_{\R^N}  e^{-\alpha (\sum_{i=1}^{N} v_{i})^2-\frac{1}{2(1+c)}\sum_{i=1}^{N} v_{i}^2}  \underset{2 \leq i<j}{\Pi} |v_j-v_i| d v_1 \cdots d v_N
\end{equation*}
By using equation (\ref{eq:integral2}) in the appendix, this leads finally to the following equivalent as $\epsilon \to 0$
\begin{align}\label{eq:Exprecepsilon}
c_\epsilon  \underset{\epsilon \to 0}{\sim}  \frac{1}{N} \frac{\Gamma(1/2)}{\Gamma(N/2)}(1+c)^{(N-1)/2} \sigma_{\epsilon}^{N-1} e^{\frac{\sigma_{\epsilon}^2}{2}}.
\end{align}

\subsection{Computation of the moment of order 2}

In order to study the convergence, for $\epsilon \to 0$, of the Gaussian chaos $M^\epsilon(A)$ (Eq. \ref{eq:GaussChaos}), we need to consider first the second-order moment 
$
E(M^{\epsilon}(A)^2)=\frac{1}{c_{\epsilon}^2}\int_{A\times A} E (e^{X^{\epsilon}(x)}e^{X^{\epsilon}(y)}) dxdy,
$
that involves the following quantity
\begin{equation}\label{eq:EeXxeXy}
 E (e^{X^{\epsilon}(x)}e^{X^{\epsilon}(y)}) = \frac{1}{N}E\left[{\rm tr}(e^{X^{\epsilon}(x)}e^{X^{\epsilon}(y)})\right]I_N.
\end{equation}
We will show that $E(M^{\epsilon}(A)^2)$ converges to a limit as $\epsilon \to 0$. Similarly, one can also prove that the sequence $(M^{\epsilon}(A))_{\epsilon >0}$ is a $L^2$ Cauchy sequence. Indeed, if $\epsilon, \epsilon'>0$ are two positive real numbers, we can write:
\begin{equation*} 
E[ {\rm tr} \left ( (M^{\epsilon}(A)-M^{\epsilon'}(A))^2 \right ) ]=E[ {\rm tr} (M^{\epsilon}(A))^2]+E[ {\rm tr} (M^{\epsilon'}(A))^2]-2 E[ {\rm tr} (M^{\epsilon}(A)M^{\epsilon'}(A))]. 
\end{equation*}
We can then conclude along the same lines as below that $E[ {\rm tr} (M^{\epsilon}(A))^2]$ and $E[ {\rm tr} (M^{\epsilon}(A)M^{\epsilon'}(A))]$ converge to the same limit as $\epsilon, \epsilon'$ go to $0$. 

Again, the proportionality to the identity matrix in (\ref{eq:EeXxeXy})  comes from the isotropic character of matrices and we will see moreover that, because the so-defined field of matrices is homogeneous, the former quantity will depend only on $|x-y|$. The purpose of this section is to compute this quantity. We will restrict to the case $|y-x|> \epsilon$ as the case $|y-x|\leq \epsilon$, once integrated, leads to vanishing terms in the limit $\epsilon \to 0$. It requires first the derivation of the joint density of the two matrices $X^{\epsilon}(x)$ and $X^{\epsilon}(y)$. We will see indeed that the quantity will depend only on $|x-y|$. We will also notice that, contrary to the one-point density (Eq. \ref{eq:1pointDensity}) from which it can be shown that eigenvectors and eigenvalues are independent, eigenvalues at point $x$ are not only correlated to eigenvalues at point $y$, but also with eigenvectors at point $y$. This intricate correlation structure is reminiscent of the noncommutative nature of this field of matrices and is encoded in the so-called Harish-Chandra--Itzykson-Zuber integral over the orthogonal group, or angular-matrix integral, and its related moments. This is an active field of research in random matrix theory and up to now, no explicit formula are known in dimension $N\geq 3$ (see for instance \cite{BreHik03,BerEyn09,ColGui09} and references therein). Nonetheless, we will succeed to get an explicit result in the asymptotic limit $\epsilon \to 0$.

\subsubsection{Joint density of two isotropic matrices}

We consider here two isotropic symmetric matrices $X^{\epsilon}(x)=(X^{\epsilon}_{i,j}(x))_{1 \leq i,j \leq N}$ and $X^{\epsilon}(y)=(X^{\epsilon}_{i,j}(y))_{1 \leq i,j \leq N}$ as defined in section \ref{eremat}. We recall that matrix components are logarithmically correlated over space. We note $x_{i,j}=X^{\epsilon}_{i,j}(x)$ and $y_{i,j}=X^{\epsilon}_{i,j}(y)$, and in matrix form $X=X^{\epsilon}(x)$ and $Y=X^{\epsilon}(y)$.

Let us first consider the diagonal terms 
\begin{equation*}
(x_{1,1}, \cdots , x_{N,N}, y_{1,1}, \cdots , y_{N,N}).
\end{equation*} 
The covariance structure $K_{2N}$ of these elements is given by
 \begin{equation*}
 K_{2N}=\left(
\begin{array}{cc}
\sigma_{\epsilon}^2 A_N  & \sigma_{|y-x|}^2 A_N   \\
 \sigma_{|y-x|}^2 A_N  &  \sigma_{\epsilon}^2 A_N \\     
\end{array}
\right),
\end{equation*}
where $A_N= (1+c) I_{N}-cP_N$ and we recall that $\sigma_{\epsilon}^2 = \gamma^2 (\ln \frac{L}{\epsilon}+1) $ and $\sigma_{|x-y|}^2 = \gamma^2 \ln \frac{L}{|x-y|}$. We know that the inverse of $K_{2N}$ is given by
\begin{equation*}
 K_{2N}^{-1}=\frac{1}{\sigma_{\epsilon}^4-\sigma_{|y-x|}^4 }\left(
\begin{array}{cc}
\sigma_{\epsilon}^2 A_N^{-1}  & -\sigma_{|y-x|}^2 A_N^{-1}    \\
 -\sigma_{|y-x|}^2 A_N^{-1}   &  \sigma_{\epsilon}^2 A_N^{-1}  \\     
\end{array}
\right),
\end{equation*}
where $A_N^{-1}=\frac{1}{(1+c)}I_{N}+2\alpha P_N$ with $\alpha=\frac{c}{2(1+c)}\frac{1}{(1+c(1-N))}$ which leads to the following density
\begin{align*}
&f((x_{i,i})_{1\leq i\leq N};(y_{j,j})_{1\leq j\leq N})=\\
&c_N  e^{-\frac{\sigma_{\epsilon}^2/(1+c) \sum_i x_{i,i}^2 + 2 \alpha \sigma_{\epsilon}^2 (\sum_i x_{i,i})^2+\sigma_{\epsilon}^2/(1+c) \sum_i y_{i,i}^2 + 2 \alpha \sigma_{\epsilon}^2 (\sum_i y_{i,i})^2-2 \sigma_{|y-x|}^2/(1+c) \sum_i x_{i,i} y_{i,i} -4 \sigma_{|y-x|}^2 \alpha (\sum_i x_{i,i})( \sum_i y_{i,i})  }{2(\sigma_{\epsilon}^4-\sigma_{|y-x|}^4)}}
\end{align*}
where $c_N=\frac{1}{(2 \pi)^N\sqrt{\text{det}(K_{2N})}}$. Now, $\text{det}(K_{2N})=(\sigma_{\epsilon}^4-\sigma_{|y-x|}^4)^{N}(1+c)^{2(N-1)}(1+c(1-N))^2$ and therefore $c_N=\frac{1}{(2 \pi)^N (\sigma_{\epsilon}^4-\sigma_{|y-x|}^4)^{N/2}(1+c)^{(N-1)}(1+c(1-N))}$. A similar procedure can be performed for the remaining $N(N-1)$ off-diagonal terms of the two matrices. The density of the couple $(X=X^{\epsilon}(x), Y=X^{\epsilon}(y))$ is thus given by, in matrix form
\begin{align} 
& f(X,Y) =\notag \\
&\bar{c}_N   e^{- \frac{\sigma_{\epsilon}^2}{2(1+c) (\sigma_{\epsilon}^4-\sigma_{|y-x|}^4) } ({\rm tr} X^2+{\rm tr} Y^2 ) - \frac{\alpha \sigma_{\epsilon}^2}{(\sigma_{\epsilon}^4-\sigma_{|y-x|}^4) } ( ({\rm tr} X)^2+({\rm tr} Y)^2 )  + \frac{\sigma_{|y-x|}^2}{(1+c) (\sigma_{\epsilon}^4-\sigma_{|y-x|}^4) } {\rm tr} XY   +\frac{2 \alpha \sigma_{|y-x|}^2}{\sigma_{\epsilon}^4-\sigma_{|y-x|}^4}  {\rm tr} X {\rm tr} Y }\label{eq:2pointjointdensXY}
\end{align}
where $\bar{c}_N=c_N \frac{1}{ \pi ^{N(N-1)/2}(1+c)^{N(N-1)/2} (\sigma_{\epsilon}^4-\sigma_{|y-x|}^4)^{N(N-1)/4} }$. We can see in the expression of the joint density of the two matrices $X$ and $Y$ (Eq. \ref{eq:2pointjointdensXY}) two different contributions. The first one, with terms of the form ${\rm tr} X^2+{\rm tr} Y^2$ and $({\rm tr} X)^2+({\rm tr} Y)^2$, relates the density of two symmetric isotropic matrices as if they were independent. The second contribution relates an interaction term coming from the logarithmic correlation of the components. Indeed, the former vanishes if the matrices are independent, i.e. $\sigma_{|y-x|}^2=0$.

At this stage, it is convenient to introduce two i.i.d. random matrices $M=(M_{i,j})$ and $M'=(M'_{i,j})$. These random matrices are taken to be living in the Gaussian Orthogonal Ensemble (GOE), namely they are symmetric and isotropic with independent components with the following distribution: the components $(M_{i,j})_{i \leq j}$ are independent centered Gaussian variables with the following variances
\begin{equation*} 
E[  M_{i,j}^2 ]=  \frac{(1+c) (\sigma_{\epsilon}^4-\sigma_{|y-x|}^4)}{2 \sigma_{\epsilon}^2}  , \; i<j; \quad E[  M_{i,i}^2 ]=   \frac{(1+c) (\sigma_{\epsilon}^4-\sigma_{|y-x|}^4)}{ \sigma_{\epsilon}^2}. 
\end{equation*}
With this, we get the following expression for $E[F(X(x), X(y))]$, where $F$ is any functional of the two matrices $X(x)$ and $X(y)$
\begin{align*} 
& E[F(X(x), X(y))]\\ &=\frac{1}{Z} E\left[    F(M,M') e^{ - \frac{\alpha \sigma_{\epsilon}^2}{ (\sigma_{\epsilon}^4-\sigma_{|y-x|}^4) } ( ({\rm tr} M)^2+({\rm tr} M')^2 ) + \frac{\sigma_{|y-x|}^2}{(1+c) (\sigma_{\epsilon}^4-\sigma_{|y-x|}^4) } {\rm tr} MM'   +\frac{2 \alpha \sigma_{|y-x|}^2}{\sigma_{\epsilon}^4-\sigma_{|y-x|}^4}  {\rm tr} M {\rm tr} M' }   \right] ,
\end{align*}
where 
\begin{equation} \label{eq:Z2points}
Z=E\left[  e^{ - \frac{\alpha \sigma_{\epsilon}^2}{ (\sigma_{\epsilon}^4-\sigma_{|y-x|}^4) } ( ({\rm tr} M)^2+({\rm tr} M')^2 ) + \frac{\sigma_{|y-x|}^2}{(1+c) (\sigma_{\epsilon}^4-\sigma_{|y-x|}^4) } {\rm tr} MM'   +\frac{2 \alpha \sigma_{|y-x|}^2}{\sigma_{\epsilon}^4-\sigma_{|y-x|}^4}  {\rm tr} M {\rm tr} M' } \right].   
\end{equation}
By using classical theorems about isotropic matrices (see \cite{Anderson}), we know that $M=OD(\lambda) {}^tO$, $M'=O'D(\lambda') {}^tO'$  where $O$ (resp. $O'$) is uniformly distributed on the orthogonal  group of $\R^N$ and is independent of the diagonal matrix $D(\lambda)$ (resp. $D(\lambda')$) the diagonal entries of which are the eigenvalues of $M$ (resp. $M'$). 

\subsubsection{Joint density of eigenvalues of two correlated isotropic matrices}
We are interested here in computing the renormalization constant $Z$ (Eq. \ref{eq:Z2points}). To do so, we diagonalize the matrices $M$ and $M'$, and perform an integration over the remaining degrees of freedom left by the eigenvectors (see \cite{BerEyn09} for instance). We define the eigenvalues of $M$ as $\lambda=(\lambda_1,\dots,\lambda_N)\in\R^d$ and we denote the Vandermonde determinant by $\triangle(\lambda)=\prod_{1\leq i<j\leq N}|\lambda_i-\lambda_j|$. We get
\begin{align*} 
Z & =\frac{1}{R_N^{\epsilon}} \int_{\R^N \times \R^N} |\Delta(\lambda)| |\Delta(\lambda')| e^{- \frac { \sigma_{\epsilon}^2} {2 (1+c) (\sigma_{\epsilon}^4-\sigma_{|y-x|}^4)} \sum_{i=1}^N \lambda_i^2 - \frac { \sigma_{\epsilon}^2} {2 (1+c) (\sigma_{\epsilon}^4-\sigma_{|y-x|}^4)} \sum_{i=1}^N \lambda_i'^2  }   \\
& \times e^{ - \frac{\alpha \sigma_{\epsilon}^2}{ (\sigma_{\epsilon}^4-\sigma_{|y-x|}^4) } ( (\sum_{i=1}^N \lambda_i)^2+(\sum_{i=1}^N \lambda_i' )^2 )  +\frac{2 \alpha \sigma_{|y-x|}^2}{\sigma_{\epsilon}^4-\sigma_{|y-x|}^4}(\sum_{i=1}^N \lambda_i)(\sum_{i=1}^N \lambda_i' ) } J( D(\lambda) , D(\lambda')  ) d \lambda d \lambda',
\end{align*}
where $R^\epsilon_N$ is a renormalization constant such that  
\begin{align*} 
\frac{1}{R_N^{\epsilon}} \int_{\R^N \times \R^N} |\Delta(\lambda)| |\Delta(\lambda')| e^{- \frac { \sigma_{\epsilon}^2} {2 (1+c) (\sigma_{\epsilon}^4-\sigma_{|y-x|}^4)} \sum_{i=1}^N \lambda_i^2 - \frac { \sigma_{\epsilon}^2} {2 (1+c) (\sigma_{\epsilon}^4-\sigma_{|y-x|}^4)} \sum_{i=1}^N \lambda_i'^2  }   d \lambda d \lambda'=1,
\end{align*}
and $J$ is the following Harish-Chandra-Itzykson-Zuber integral \cite{BreHik03,BerEyn09,ColGui09}, also called matrix angular integral ($dO$ stands for the Haar measure on $O_N(\R)$)
\begin{equation*}
J( D(\lambda) , D(\lambda')  )= \int_{O_N(\R)}   e^{ \frac{\sigma_{|y-x|}^2}{(1+c) (\sigma_{\epsilon}^4-\sigma_{|y-x|}^4) } {\rm tr} \:  D(\lambda) O  D(\lambda') O^{-1}  }  dO,
\end{equation*}
obtained while integrating over the eigenvectors that enter in the term ${\rm tr} MM'$ of Eq. \ref{eq:Z2points}.  We make the change of variables $u_i=\frac{\sigma_{\epsilon}}{\sqrt{\sigma_{\epsilon}^4-\sigma_{|y-x|}^4}} \lambda_i$, $u_i'=\frac{\sigma_{\epsilon}}{\sqrt{\sigma_{\epsilon}^4-\sigma_{|y-x|}^4}} \lambda_i'$ (set $\gamma_{\epsilon}=\frac{\sqrt{\sigma_{\epsilon}^4-\sigma_{|y-x|}^4}} {\sigma_{\epsilon}}$) and get:
\begin{align*} 
Z & =\frac{\gamma_{\epsilon}^{N(N+1)}}{R_N^{\epsilon}} \int_{\R^N \times \R^N} |\Delta(u)| |\Delta(u')| e^{- \frac {1} {2 (1+c)} \sum_{i=1}^N u_i^2 - \frac { 1} {2 (1+c)} \sum_{i=1}^N u_i'^2  }   \\
& \times e^{ - \alpha ( (\sum_{i=1}^N u_i)^2+(\sum_{i=1}^N u_i' )^2 )  +\frac{2 \alpha \sigma_{|y-x|}^2}{\sigma_{\epsilon}^2}(\sum_{i=1}^N u_i)(\sum_{i=1}^N u_i' ) } \overline{J}(   D(u) ,   D(u')  ) d u d u',
\end{align*}
where we have set
\begin{equation*}
 \overline{J}(  D(u) ,  D(u')  )  = \int_{O_N(\R)}   e^{  \frac{1}{1+c} \frac{ \sigma_{|y-x|}^2}{\sigma_{\epsilon}^2} \sum_{i,j=1}^N u_i u_j' |O_{i,j}|^2 }  dO.
\end{equation*}
Therefore, since $\overline{J}(  D(u) ,  D(u')  )$ converges pointwise towards $1$ as $\epsilon\to 0$, we can use the Lebesgue dominated convergence theorem to get the following equivalent as $\epsilon \to 0$
\begin{align*} 
Z &  \underset{\epsilon \to 0}{\sim} \frac{\gamma_{\epsilon}^{N(N+1)}}{R_N^{\epsilon}} \int_{\R^N \times \R^N} |\Delta(u)| |\Delta(u')| e^{- \frac {1} {2 (1+c)} \sum_{i=1}^N u_i^2 - \frac { 1} {2 (1+c)} \sum_{i=1}^N u_i'^2  }   \\
& \times e^{ - \alpha ( (\sum_{i=1}^N u_i)^2+(\sum_{i=1}^N u_i' )^2 ) }  d u d u',
\end{align*}
that is straightforward to compute (see the appendix).

\subsubsection{Two-points correlation structure of the matrix chaos}

We want to get an equivalent as $\epsilon\to 0$ of the quantity given in Eq. \ref{eq:EeXxeXy}. To do so, we consider the following quantity
\begin{equation*}
\bar{Z}= E[ {\rm tr} (e^{M} e^{M'}) e^{ - \frac{\alpha \sigma_{\epsilon}^2}{ (\sigma_{\epsilon}^4-\sigma_{|y-x|}^4) } ( ({\rm tr} M)^2+({\rm tr} M')^2 ) + \frac{\sigma_{|y-x|}^2}{(1+c) (\sigma_{\epsilon}^4-\sigma_{|y-x|}^4) } {\rm tr} MM'   +\frac{2 \alpha \sigma_{|y-x|}^2}{\sigma_{\epsilon}^4-\sigma_{|y-x|}^4}  {\rm tr} M {\rm tr} M' } ].
\end{equation*}
In the same spirit as formerly, we diagonalize the matrices $M$ and $M'$ and perform the integration over the eigenvectors. We get
\begin{align*}
 \bar{Z}   & =    \frac{1}{R_N^{\epsilon}} \int_{\R^N \times \R^N} |\Delta(\lambda)| |\Delta(\lambda')| e^{- \frac { \sigma_{\epsilon}^2} {2 (1+c) (\sigma_{\epsilon}^4-\sigma_{|y-x|}^4)} \sum_{i=1}^N \lambda_i^2 - \frac { \sigma_{\epsilon}^2} {2 (1+c) (\sigma_{\epsilon}^4-\sigma_{|y-x|}^4)} \sum_{i=1}^N \lambda_i'^2  }   \\
& \times e^{ - \frac{\alpha \sigma_{\epsilon}^2}{ (\sigma_{\epsilon}^4-\sigma_{|y-x|}^4) } ( (\sum_{i=1}^N \lambda_i)^2+(\sum_{i=1}^N \lambda_i' )^2 )  +\frac{2 \alpha \sigma_{|y-x|}^2}{\sigma_{\epsilon}^4-\sigma_{|y-x|}^4}(\sum_{i=1}^N \lambda_i)(\sum_{i=1}^N \lambda_i' ) } I( D(\lambda) , D(\lambda')  ) d \lambda d \lambda',
\end{align*}
where $I$ is the following moment of the angular integral
\begin{equation*}
I ( D(\lambda) , D(\lambda')  )= \int_{O_N(\R)}  {\rm tr}( e^{D(\lambda) } O  e^{D(\lambda') }  O^{-1} ) e^{ \frac{\sigma_{|y-x|}^2}{(1+c) (\sigma_{\epsilon}^4-\sigma_{|y-x|}^4) } {\rm tr} \:  D(\lambda) O  D(\lambda') O^{-1}  }  dO.
\end{equation*}
We make again the change of variables $u_i=\frac{\sigma_{\epsilon}}{\sqrt{\sigma_{\epsilon}^4-\sigma_{|y-x|}^4}} \lambda_i$, $u_i'=\frac{\sigma_{\epsilon}}{\sqrt{\sigma_{\epsilon}^4-\sigma_{|y-x|}^4}} \lambda_i'$ (set $\gamma_{\epsilon}=\frac{\sqrt{\sigma_{\epsilon}^4-\sigma_{|y-x|}^4}} {\sigma_{\epsilon}}$)
\begin{align*} 
\bar{Z} & = \sum_{i,j=1}^N \frac{\gamma_{\epsilon}^{N(N+1)}}{R_N^{\epsilon}} \int_{\R^N \times \R^N} |\Delta(u)| |\Delta(u')| e^{- \frac {1} {2 (1+c)} \sum_{k=1}^N u_k^2 - \frac { 1} {2 (1+c)} \sum_{k=1}^N u_k'^2  }   \\
& \times e^{ - \alpha( (\sum_{k=1}^N u_k)^2+(\sum_{k=1}^N u_k' )^2 )  +\frac{2 \alpha \sigma_{|y-x|}^2}{\sigma_{\epsilon}^2}(\sum_{k=1}^N u_k)(\sum_{k=1}^N u_k' ) } e^{\gamma_{\epsilon}(u_i+u_j')} \overline{I}_{i,j}(  D(u) ,   D(u')  ) d u d u',
\end{align*}
where we have set
\begin{equation*}
 \overline{I}_{i,j}( D(u) , D(u')  )  = \int_{O_N(\R)} |O_{i,j}|^2   e^{  \frac{1}{1+c} \frac{ \sigma_{|y-x|}^2}{\sigma_{\epsilon}^2} \sum_{k,k'=1}^N u_k u_{k'}' |O_{k,k'}|^2 }  dO,
\end{equation*}
known as the Morozov moment \cite{BerEyn09}. We make the following change of variables in the above integral: $u_i=v_i+\gamma_{\epsilon},\; u_k=v_k-c\gamma_{\epsilon}$ for $k \not = i $ and $u_j'=v_j'+\gamma_{\epsilon},\;  u_k'=v_k'-c\gamma_{\epsilon}$ for $k \not = j$. We obtain the following equivalent
\begin{align*} 
\bar{Z} & \underset{\epsilon \to 0}{\sim} \sum_{i,j=1}^N \frac{\gamma_{\epsilon}^{N(N+1)}e^{\sigma_{\epsilon}^2}(1+c)^{2(N-1)} \sigma_{\epsilon}^{2(N-1)}}{R_N^{\epsilon}}   I_{i,j} \int_{\R^N \times \R^N} |\Delta_i(v)| |\Delta_j(v')|\\ 
&\times e^{- \frac {1} {2 (1+c)} \sum_{k=1}^N v_k^2 - \frac { 1} {2 (1+c)} \sum_{k=1}^N v_k'^2  - \alpha ( (\sum_{k=1}^N v_k)^2+(\sum_{k=1}^N v_k' )^2 )  + 2 \alpha \sigma_{|y-x|}^2(1+c(1-N))^2 }  d v d v',
\end{align*}
where $|\Delta_i(v)|= \prod_{l<l', l,l' \not = i } | v_l-v_{l'} |  $ and:
\begin{align*} 
I_{i,j} & = \int_{O_n(\R)} |O_{i,j}|^2   e^{  \frac{1}{1+c}  \sigma_{|y-x|}^2  \sum_{k,k'=1}^N (-c+(1+c)1_{k=i}) (-c+(1+c)1_{k'=j}) |O_{k,k'}|^2 }  dO   \\
& =  e^{  \sigma_{|y-x|}^2  (   \frac{c^2 N}{1+c} -2c)}   \int_{O_N(\R)} |O_{1,1}|^2   e^{  \sigma_{|y-x|}^2(1+c) |O_{1,1}|^2  }  dO,   \\ 
\end{align*}
which is independent of $i,j$. Therefore, we get
\begin{align*} 
\bar{Z} & \underset{\epsilon \to 0}{\sim} N^2 \frac{\gamma_{\epsilon}^{N(N+1)}e^{\sigma_{\epsilon}^2}(1+c)^{2(N-1)} \sigma_{\epsilon}^{2(N-1)}}{R_N^{\epsilon}}   I_{1,1} \int_{\R^N \times \R^N} |\Delta_1(v)| |\Delta_1(v')| \\
&e^{- \frac {1} {2 (1+c)} \sum_{k=1}^N v_k^2 - \frac { 1} {2 (1+c)} \sum_{k=1}^N v_k'^2  - \alpha ( (\sum_{k=1}^N v_k)^2+(\sum_{k=1}^N v_k' )^2 )  + 2 \alpha \sigma_{|y-x|}^2(1+c(1-N))^2 }  d v d v'.
\end{align*}
In conclusion, we get
\begin{equation*} 
\bar{Z} / Z \underset{\epsilon \to 0}{\sim}    (1+c)^{N-1} (\frac{\Gamma(1/2)}{\Gamma(N/2)})^2  e^{\sigma_{\epsilon}^2} \sigma_{\epsilon}^{2(N-1)}    e^{-c \sigma_{|y-x|}^2}  \int_{O_N(\R)} |O_{1,1}|^2   e^{  \sigma_{|y-x|}^2(1+c) |O_{1,1}|^2  }  dO.
\end{equation*}
Including furthermore the normalization constant $c_{\epsilon}$ (Eq. \ref{eq:Exprecepsilon}), we get
\begin{equation*} 
\bar{Z} / (Z c_{\epsilon}^2) \underset{\epsilon \to 0}{\sim}    N^2  e^{-c \sigma_{|y-x|}^2}  \int_{O_N(\R)} |O_{1,1}|^2   e^{  \sigma_{|y-x|}^2(1+c) |O_{1,1}|^2  }  dO. 
\end{equation*}

\subsubsection{Computation of the moment of order $2$}
From the above subsections, we deduce that
\begin{equation*}
E(tr M^{\epsilon}(A)^2) \underset{\epsilon \to 0}{\rightarrow} N^2 \int_{A\times A} e^{-c \sigma_{|y-x|}^2}  \int_{O_N(\R)} |O_{1,1}|^2   e^{  \sigma_{|y-x|}^2(1+c) |O_{1,1}|^2  }  dO \:  dxdy
\end{equation*}

We recall that the law of $|O_{1,1}|^2$ is the one of the square of one component of a vector uniformly distributed on the unit sphere, and has thus a density given by (see Lemma \ref{haargauss})
\begin{equation*}
f(v)=   \frac{\Gamma(N/2)}{\Gamma(1/2) \Gamma((N-1)/2)}v^{-1/2} (1-v)^{(N-3)/2} .
\end{equation*}
We get the following equivalent as $|y-x| \to 0$:
\begin{equation*}
N^2 e^{-c \sigma_{|y-x|}^2}  \int_{O_N(\R)} |O_{1,1}|^2   e^{  \sigma_{|y-x|}^2(1+c) |O_{1,1}|^2  }  dO \underset{|y-x| \to 0}{\sim}   N^2 \frac{\Gamma(N/2)}{\Gamma(1/2)} \frac{e^{ \sigma_{|y-x|}^2} }{(1+c)^{(N-1)/2}  \sigma_{|y-x|}^{N-1}  },
\end{equation*}
which entails (\ref{firstequivalent}).

\subsection{Computation of the moment of order $k$}

We are interested here in studying the convergence, when $\epsilon \to 0$, of the Gaussian chaos $M^\epsilon(A)$ (Eq. \ref{eq:GaussChaos}) for higher order moments such as, $k\in \mathbb N$,
\begin{equation*}
E(M^{\epsilon}(A))^k=\frac{1}{c_{\epsilon}^k}\int_{A^k} E \left(\prod_{1\leq i\leq k}e^{X^{\epsilon}(x_i)}\right) dx_1\cdots dx_k,
\end{equation*}
that involves the following quantity
\begin{equation}\label{eq:EeXxi}
E \left(\prod_{1\leq i\leq k}e^{X^{\epsilon}(x_i)}\right) = \frac{1}{N} E\left[{\rm tr}\prod_{1\leq i\leq k}e^{X^{\epsilon}(x_i)}\right]I_N.
\end{equation}
In this subsection, we will suppose that $k$ is an integer greater or equal to $2$ such that $k<\frac{2d}{\gamma^2}$. This condition ensures that all the integrals we consider below are finite and that one can apply the dominated convergence theorem to justify the inversions between limit and integral we will perform with no further justification. To generalize the former calculations in the case $k=2$, we will first derive the joint density of $k$-matrices $(X^{\epsilon}(x_i))_{1\leq i\leq k}$. A generalized version  to $k$-points of the Harish-Chandra-Itzykson-Zuber integral enters the expression of the density. An exact evaluation of these integrals remains an open issue. As far as we know, only their behavior in the asymptotic  limit of large matrices ($N\rightarrow +\infty$) has been considered in the literature \cite{ColGui09}. Nonetheless, a logarithmic equivalent of the quantity of interest (Eq. \ref{eq:EeXxi}) can be obtained and allows us to show the multifractal behavior of the multiplicative chaos (i.e. $\zeta(k)$ is a non linear function of the order $k$, see theorem \ref{Theo:Ncase}).

\subsubsection{Joint density of $k$ isotropic Gaussian matrices}

We consider here $k$ isotropic Gaussian matrices $\left( X^\epsilon(x_i)\right)_{1\leq i\leq k}$. The ensemble made of the $kN$ diagonal terms, i.e. 
\begin{equation*}
(X^\epsilon_{1,1}(x_1), \cdots , X^\epsilon_{N,N}(x_1), \cdots, X^\epsilon_{1,1}(x_k), \cdots , X^\epsilon_{N,N}(x_k)),
\end{equation*} 
has covariance structure $K_{kN}$
 \begin{equation*}
 K_{kN}=\left(
\begin{array}{ccccc}
\sigma_{\epsilon}^2 A_N &  \sigma_{|x_1-x_2|}^2  A_N &  \cdots &  \cdots  & \sigma_{|x_1-x_k|}^2 A_N   \\
 \sigma_{|x_2-x_1|}^2 A_N  & \sigma_{\epsilon}^2 A_N & \cdots &   \cdots  &  \sigma_{|x_2-x_k|}^2 A_N \\     
   \cdots  &  \cdots  & \cdots & \cdots  & \cdots  \\
 \sigma_{|x_{k-1}-x_1|}^2 A_N    &   \sigma_{|x_{k-1}-x_2|}^2 A_N   & \cdots &   \sigma_{\epsilon}^2 A_N &  \sigma_{|x_{k-1}-x_k|}^2 A_N  \\
 \sigma_{|x_k-x_1|}^2 A_N  &    \sigma_{|x_k-x_2|}^2 A_N & \cdots &  \sigma_{|x_{k}-x_{k-1}|}^2 A_N  & \sigma_{\epsilon}^2 A_N \\
\end{array}
\right),
\end{equation*}
where again, $A_N= (1+c) I_{N}-cP_N$. We know that the inverse of $K_{kN}$ is approximately given by ($\epsilon \to 0$)
\begin{equation*}
 K_{kN}^{-1}=\frac{1}{\sigma_{\epsilon}^4 }\left(
\begin{array}{ccccc}
\sigma_{\epsilon}^2 A_N^{-1} &  -\sigma_{|x_1-x_2|}^2  A_N^{-1} &  \cdots &  \cdots  & - \sigma_{|x_1-x_k|}^2 A_N^{-1}   \\
- \sigma_{|x_2-x_1|}^2 A_N^{-1}  & \sigma_{\epsilon}^2 A_N^{-1} & \cdots &   \cdots  & - \sigma_{|x_2-x_k|}^2 A_N^{-1} \\     
   \cdots  &  \cdots  & \cdots & \cdots  & \cdots  \\
 - \sigma_{|x_{k-1}-x_1|}^2 A_N^{-1}    & -  \sigma_{|x_{k-1}-x_2|}^2 A_N^{-1}   & \cdots &   \sigma_{\epsilon}^2 A_N^{-1} & - \sigma_{|x_{k-1}-x_k|}^2 A_N^{-1}  \\
- \sigma_{|x_k-x_1|}^2 A_N^{-1}  &   - \sigma_{|x_k-x_2|}^2 A_N^{-1} & \cdots & - \sigma_{|x_{k}-x_{k-1}|}^2 A_N^{-1}  & \sigma_{\epsilon}^2 A_N^{-1} \\
\end{array}
\right),
\end{equation*}
where $A_N^{-1}=\frac{1}{(1+c)}I_{N}+2\alpha P_N$, with $\alpha=\frac{c}{2(1+c)}\frac{1}{(1+c(1-N))}$.
The density of diagonal components, considering the $N$-dimensional vector $X^{(l)}=(X^\epsilon_{1,1}(x_l), \cdots , X^\epsilon_{N,N}(x_l))$, is thus given by
\begin{align*}
f(X^{(1)}, \cdots , X^{(k)}) =c_N e^{-\frac{1}{2 \sigma_{\epsilon}^4}  \sum_{i,j=1}^{k}  (\delta_{i,j}\sigma_{\epsilon}^2 -(1-\delta_{i,j})\sigma_{|x_{i}-x_j|}^2  )   ^{t} X^{(i)}  (\frac{1}{(1+c)}I_{N}+2 \alpha P_N) X^{(j)} }
\end{align*}
where $c_N=\frac{1}{(2 \pi)^{kN/2}\sqrt{\text{det}(K_{kN})}}$. 
For the off-diagonal terms, the situation is simpler. If $i<j$, the covariance matrix of the vector $(X^\epsilon_{i,j}(x_1), \cdots, X^\epsilon_{i,j}(x_k))$, which is independent on all the remaining diagonal and off-diagonal components, is
 \begin{equation*}
 \frac{1+c}{2}
 \left(
\begin{array}{ccccc}
\sigma_{\epsilon}^2  &  \sigma_{|x_1-x_2|}^2   &  \cdots &  \cdots  & \sigma_{|x_1-x_k|}^2    \\
 \sigma_{|x_2-x_1|}^2   & \sigma_{\epsilon}^2  & \cdots &   \cdots  &  \sigma_{|x_2-x_k|}^2  \\     
   \cdots  &  \cdots  & \cdots & \cdots  & \cdots  \\
 \sigma_{|x_{k-1}-x_1|}^2     &   \sigma_{|x_{k-1}-x_2|}^2    & \cdots &   \sigma_{\epsilon}^2  &  \sigma_{|x_{k-1}-x_k|}^2   \\
 \sigma_{|x_k-x_1|}^2   &    \sigma_{|x_k-x_2|}^2  & \cdots &  \sigma_{|x_{k}-x_{k-1}|}^2   & \sigma_{\epsilon}^2  \\
\end{array}
\right),
\end{equation*}
whose inverse is approximately given by ($\epsilon \to 0$)
\begin{equation*}
\frac{2}{(1+c)\sigma_{\epsilon}^4 }\left(
\begin{array}{ccccc}
\sigma_{\epsilon}^2  &  -\sigma_{|x_1-x_2|}^2  &  \cdots &  \cdots  & - \sigma_{|x_1-x_k|}^2    \\
- \sigma_{|x_2-x_1|}^2   & \sigma_{\epsilon}^2  & \cdots &   \cdots  & - \sigma_{|x_2-x_k|}^2  \\     
   \cdots  &  \cdots  & \cdots & \cdots  & \cdots  \\
 - \sigma_{|x_{k-1}-x_1|}^2     & -  \sigma_{|x_{k-1}-x_2|}^2   & \cdots &   \sigma_{\epsilon}^2  & - \sigma_{|x_{k-1}-x_k|}^2   \\
- \sigma_{|x_k-x_1|}^2   &   - \sigma_{|x_k-x_2|}^2  & \cdots & - \sigma_{|x_{k}-x_{k-1}|}^2   & \sigma_{\epsilon}^2  \\
\end{array}
\right).
\end{equation*}
This leads to the following density, using the notations $x_{i,j}^{(r)}=X^\epsilon_{i,j}(x_r)$
\begin{align*}
f(x_{i,j}^{(1)}, \cdots, x_{i,j}^{(k)}) = k_N  e^{-\frac{1}{(1+c ) \sigma_{\epsilon}^4}  \sum_{r,l=1}^{k}  (\delta_{r,l}\sigma_{\epsilon}^2 -(1-\delta_{r,l})\sigma_{|x_{r}-x_{l} |}^2  )  x_{i,j}^{(r)} x_{i,j}^{(l)} } .
\end{align*}
Therefore, we get the following density for the $k$ matrices (we omit superscript $\epsilon$ for the sake of clarity)
\begin{align*}
f(X(x_1), \cdots , X(x_k)) &= \bar{c}_N e^{-\frac{1}{2 \sigma_{\epsilon}^4}  \sum_{r,l=1}^{k}  (\delta_{r,l}\sigma_{\epsilon}^2 -(1-\delta_{r,l})\sigma_{|x_{r}-x_l|}^2  )   ^{t} X^{(r)}  (\frac{1}{(1+c)}I_{N}+2 \alpha P) X^{(l)} } \\
& \times e^{- \sum_{i<j} \frac{1}{(1+c ) \sigma_{\epsilon}^4}  \sum_{r,l=1}^{k}  (\delta_{r,l}\sigma_{\epsilon}^2 -(1-\delta_{r,l})\sigma_{|x_{r}-x_{l} |}^2  )  x_{i,j}^{(r)} x_{i,j}^{(l)} },
\end{align*}
which we rewrite under matrix notation
\begin{align*}
f(X(x_1), \cdots , X(x_k)) &=\bar{c}_N e^{-\frac{1}{2 (1+c) \sigma_{\epsilon}^2} \sum_{r=1}^k tr( X(x_r)^2  )-\frac{\alpha}{ \sigma_{\epsilon}^2} \sum_{r=1}^k (tr  X(x_r)  )^2  }  \\
& \times   e^{  \frac{1}{(1+c) \sigma_{\epsilon}^4}  \sum_{r<l}  \sigma_{|x_r-x_l|}^2  tr  X(x_r) X(x_l)  +\frac{2 \alpha}{\sigma_{\epsilon}^4}  \sum_{r<l}  \sigma_{|x_r-x_l|}^2  tr X(x_r)tr X(x_l)   }   .
\end{align*}

We introduce $k$ i.i.d. random matrices $M^{(l)}=(M_{i,j}^{(l)})$ pertaining to the GOE ensemble. These random matrices are symmetric and isotropic with independent components with the following distribution: the components $(M^{(l)}_{i,j})_{i \leq j}$ are independent centered Gaussian variables with the following variances
\begin{equation*} 
E[  (M^{(l)}_{i,j})^2 ]=  \frac{1+c }{2} \sigma_{\epsilon}^2  , \; i<j; \quad E[  (M^{(l)}_{i,i})^2 ]=   (1+c)  \sigma_{\epsilon}^2.
\end{equation*}
With this, we get the following expression for the expectation of any functional $F(X^\epsilon(x_1), \cdots , X^\epsilon(x_k))$ of the $k$ matrices $ X^\epsilon(x_1), \cdots , X^\epsilon(x_k) $
\begin{align*} 
&  E[F(X^\epsilon(x_1), \cdots , X^\epsilon(x_k))] \\ 
&=\frac{ E[    F(M^{(1)}, \cdots ,M^{(k)}) e^{ - \frac{\alpha}{ \sigma_{\epsilon}^2 } \sum_{r=1}^k ({\rm tr} M^{(r)})^2 + \frac{1}{(1+c) \sigma_{\epsilon}^4}  \sum_{r<l}  \sigma_{|x_r-x_l|}^2  tr  M^{(r)} M^{(l)}  + \frac{2 \alpha}{\sigma_{\epsilon}^4}  \sum_{r<l}  \sigma_{|x_r-x_l|}^2  tr M^{(r)}  tr M^{(l)} }    ] }{ Z },
\end{align*}
where: 
\begin{equation*}
Z= E[   e^{ - \frac{\alpha}{ \sigma_{\epsilon}^2 } \sum_{r=1}^k ({\rm tr} M^{(r)})^2 + \frac{1}{(1+c) \sigma_{\epsilon}^4}  \sum_{r<l}  \sigma_{|x_r-x_l|}^2  tr  M^{(r)} M^{(l)}  + \frac{2 \alpha}{\sigma_{\epsilon}^4}  \sum_{r<l}  \sigma_{|x_r-x_l|}^2  tr M^{(r)}  tr M^{(l)} }    ]. 
\end{equation*}
By using classical theorems about isotropic matrices (see \cite{Anderson}), we know that, for each $r$, $M^{(r)}=O^{(r)} D(\lambda^{(r)} ) {}^tO^{(r)}$ where $O^{(r)}$ is uniformly distributed on the orthogonal  group of $\R^N$ and is independent of the diagonal matrix $D(\lambda^{(r)} )$  the diagonal entries of which are the eigenvalues of $M^{(r)}$. 

\subsubsection{Joint density of eigenvalues of $k$ isotropic Gaussian matrices and computation of the renormalization}
We start by computing $Z$. For $\lambda=(\lambda_1,\dots,\lambda_N)\in\R^d$, we denote the Vandermonde determinant by  $\triangle(\lambda)=\prod_{1\leq i<j\leq N}|\lambda_i-\lambda_j|$. We get
\begin{align*} 
Z & =\frac{1}{R_N^{\epsilon}} \int_{\R^{kN}} \Pi_{r=1}^k |\Delta(\lambda^{(r)})|  e^{- \frac {1} {2 (1+c) \sigma_{\epsilon}^2} \sum_{r=1}^k \sum_{i=1}^N (\lambda_i^{(r)})^2 - \frac{\alpha}{ \sigma_{\epsilon}^2 } \sum_{r=1}^k  (\sum_{i=1}^N \lambda_i^{(r)})^2 }   \\
& \times e^{   \frac{2 \alpha}{\sigma_{\epsilon}^4} \sum_{r<l}\sigma_{|x_r-x_l|}^2 (\sum_{i=1}^N \lambda_i^{(r)})(\sum_{i=1}^N \lambda_i^{(l)} ) } J( D(\lambda^{(1)}), \cdots ,D(\lambda^{(k)})  ) d \lambda^{(1)} \cdots d \lambda^{(k)},
\end{align*}
where $R_N^{\epsilon}$ is a renormalization constant such that  
\begin{align*} 
\frac{1}{R_N^{\epsilon}} \int_{\R^{kN}} \Pi_{r=1}^k |\Delta(\lambda^{(r)})|  e^{- \frac {1} {2 (1+c) \sigma_{\epsilon}^2} \sum_{r=1}^k \sum_{i=1}^N (\lambda_i^{(r)})^2  } =1 ,
\end{align*}
and $J$ is the following angular integral: ($dO^{(r)}$ stands for the Haar measure on $O_N(\R)$)
\begin{align*}
&J(  D(\lambda^{(1)}), \cdots ,D(\lambda^{(k)})  )= \\
&\int_{O_N(\R)^k}   e^{ \frac{1}{(1+c)\sigma_{\epsilon}^4 }\sum_{r<l} \sigma_{|x_r-x_l|}^2 {\rm tr} \: O^{(r)} D(\lambda^{(r)}) {}^t O^{(r)} O^{(l)}  D(\lambda^{(l)}) {}^t O^{(l)}  }  dO^{(1)} \cdots dO^{(k)}.
\end{align*}
We make the  change of variables $u_i^{(r)}=\frac{\lambda_i^{(r)}}{\sigma_{\epsilon}}$
\begin{align*} 
Z & =\frac{\sigma_{\epsilon}^{N(N+1)k/2}}{R_N^{\epsilon}} \int_{\R^{kN}}  \Pi_{r=1}^k |\Delta(u^{(r)})|  e^{- \frac {1} {2 (1+c)} \sum_{r=1}^k \sum_{i=1}^N (u_i^{(r)})^2  - \alpha \sum_{r=1}^k  (\sum_{i=1}^N u_i^{(r)})^2}   \\
& \times e^{ \frac{2 \alpha}{\sigma_{\epsilon}^2} \sum_{r<l}\sigma_{|x_r-x_l|}^2 (\sum_{i=1}^N u_i^{(r)})(\sum_{i=1}^N u_i^{(l)} ) } \overline{J}( D(u^{(1)}), \cdots ,D(u^{(l)})  ) d u^{(1)} \cdots d u^{(l)},
\end{align*}
where we have set
\begin{align*}
&\overline{J}( D(u^{(1)}), \cdots ,D(u^{(k)})  ) =\\
&\int_{O_N(\R)^k}   e^{ \frac{1}{(1+c)\sigma_{\epsilon}^2 }\sum_{r<l} \sigma_{|x_r-x_l|}^2 {\rm tr} \: O^{(r)} D(u^{(r)}) {}^t O^{(r)} O^{(l)}  D(u^{(l)}) {}^t O^{(l)}  }  dO^{(1)} \cdots dO^{(k)}.
\end{align*}
Therefore, since $\overline{J}( D(u^{(1)}), \cdots ,D(u^{(k)})  )$ converges pointwise towards $1$ as $\epsilon\to 0$, we can use the Lebesgue theorem to get the following equivalent as $\epsilon \to 0$
\begin{align*} 
Z &  \underset{\epsilon \to 0}{\sim} \frac{\sigma_{\epsilon}^{N(N+1)k/2}}{R_N^{\epsilon}} \int_{\R^{kN}}  \Pi_{r=1}^k |\Delta(u^{(r)})|  e^{- \frac {1} {2 (1+c)} \sum_{r=1}^k \sum_{i=1}^N (u_i^{(r)})^2  }   \\
& \times e^{ - \alpha \sum_{r=1}^k  (\sum_{i=1}^N u_i^{(r)})^2 }d u^{(1)} \cdots d u^{(k)}.
\end{align*}

\subsubsection{$k$-point correlation structure of the multiplicative chaos}
For $i \leq j$, we want to get an equivalent as $\epsilon\to 0$ of the following quantity
\begin{equation*}
\bar{Z}= E[ (\Pi_{r=1}^k e^{M^{(r)}})_{i,j} e^{ - \frac{\alpha}{ \sigma_{\epsilon}^2 } \sum_{r=1}^k ({\rm tr} M^{(r)})^2 + \frac{1}{(1+c) \sigma_{\epsilon}^4}  \sum_{r<l}  \sigma_{|x_r-x_l|}^2  {\rm tr}   M^{(r)} M^{(l)}  + \frac{2 \alpha}{\sigma_{\epsilon}^4}  \sum_{r<l}  \sigma_{|x_r-x_l|}^2  {\rm tr} M^{(r)}  {\rm tr}  M^{(l)}} ].
\end{equation*}
We get
\begin{align*}
 \bar{Z}   & =\frac{1}{R_N^{\epsilon}} \int_{\R^{kN}} \Pi_{r=1}^k |\Delta(\lambda^{(r)})|  e^{- \frac {1} {2 (1+c) \sigma_{\epsilon}^2} \sum_{r=1}^k \sum_{i=1}^N (\lambda_i^{(r)})^2  - \frac{\alpha}{ \sigma_{\epsilon}^2 } \sum_{r=1}^k  (\sum_{i=1}^N \lambda_i^{(r)})^2 }   \\
& \times e^{ \frac{2 \alpha}{\sigma_{\epsilon}^4} \sum_{r<l}\sigma_{|x_r-x_l|}^2 (\sum_{i=1}^N \lambda_i^{(r)})(\sum_{i=1}^N \lambda_i^{(l)} ) } I( D(\lambda^{(1)}), \cdots ,D(\lambda^{(k)})  ) d \lambda^{(1)} \cdots d \lambda^{(k)},
\end{align*}
where $I$ is the following angular integral
\begin{align*}
& I(  D(\lambda^{(1)}), \cdots ,D(\lambda^{(k)})  )   = \int_{O_N(\R)^k} (\Pi_{r=1}^k O^{(r)} e^{D(\lambda^{(r)})} {}^t O^{(r)})_{i,j}   \\
&\times e^{ \frac{1}{(1+c)\sigma_{\epsilon}^4 }\sum_{r<l} \sigma_{|x_r-x_l|}^2 {\rm tr} \: O^{(r)} D(\lambda^{(r)}) {}^t O^{(r)} O^{(l)}  D(\lambda^{(l)}) {}^t O^{(l)}  }  dO^{(1)} \cdots dO^{(k)}.
\end{align*}
We make the following change of variables $u_i^{(r)}=\frac{ \lambda_i^{(r)}}{\sigma_{\epsilon}}$
\begin{align*}
 \bar{Z}   & = \sum_{j_1, \cdots, j_k=1}^{N} \frac{\sigma_{\epsilon}^{N(N+1)k/2}}{R_N^{\epsilon}} \int_{\R^{kN}} \Pi_{r=1}^k |\Delta(u^{(r)})|  e^{- \frac {1} {2 (1+c)} \sum_{r=1}^k \sum_{i=1}^N (u_i^{(r)})^2  }   \\
& \times e^{ - \alpha \sum_{r=1}^k  (\sum_{i=1}^N u_i^{(r)})^2  +\frac{2 \alpha}{\sigma_{\epsilon}^2} \sum_{r<l}\sigma_{|x_r-x_l|}^2 (\sum_{i=1}^N u_i^{(r)})(\sum_{i=1}^N u_i^{(l)} ) } e^{\sigma_{\epsilon}(u_{j_1}^{(1)}+\cdots + u_{j_k}^{(k)})} \\
&  \times \sum_{\underset{l_0=i; l_k=j}{l_1, \cdots, l_{k-1}=1}}^N I_{l_0, l_1, \cdots, l_{k-1},l_{k}}^{j_1, \cdots, j_k}( D(u^{(1)}), \cdots ,D(u^{(k)})  ) d u^{(1)} \cdots d u^{(k)},
\end{align*}
where we have set
\begin{align*}
& I_{l_0, l_1, \cdots, l_{k-1},l_{k}}^{j_1, \cdots, j_k}( D(u^{(1)}), \cdots ,D(u^{(k)})  ) = \int_{O_N(\R)^k} (\Pi_{r=1}^k O_{l_{r-1},j_r}^{(r)} O_{l_{r},j_r}^{(r)} )  \\
&\times  e^{ \frac{1}{(1+c)\sigma_{\epsilon}^2 }\sum_{r<l} \sigma_{|x_r-x_l|}^2 {\rm tr} \: O^{(r)} D(u^{(r)}) {}^t O^{(r)} O^{(l)}  D(u^{(l)}) {}^t O^{(l)}  }  dO^{(1)} \cdots dO^{(k)} .
\end{align*}

We make the following change of variables in the above integral for $1 \leq r \leq k$: $u_{j_r}^{(r)}=v_{j_r}^{(r)}+\sigma_{\epsilon},\; u_k^{(r)}=v_k^{(r)}-c\sigma_{\epsilon} \: k \not = j_r $. We obtain the following equivalent
\begin{align*}
 \bar{Z}   & \underset{\epsilon \to 0}{\sim}  \sum_{j_1, \cdots, j_k=1}^{N} \frac{\sigma_{\epsilon}^{N(N+1)k/2} e^{k\sigma_{\epsilon}^2/2} (1+c)^{(N-1)k}\sigma_{\epsilon}^{(N-1)k}}{R_N^{\epsilon}} \int_{\R^{kN}} \Pi_{r=1}^k |\Delta_{j_r}(v^{(r)})|    \\
& \times e^{- \frac {1} {2 (1+c)} \sum_{r=1}^k \sum_{i=1}^N (v_i^{(r)})^2 - \alpha \sum_{r=1}^k  (\sum_{i=1}^N v_i^{(r)})^2  +2 \alpha (1+c(1-N))^2 \sum_{r<l}\sigma_{|x_r-x_l|}^2  } \\
&  \times \sum_{\underset{l_0=i; l_k=j}{l_1, \cdots, l_{k-1}=1}}^N \bar{I}_{l_0, l_1, \cdots, l_{k-1},l_{k}}^{j_1, \cdots, j_k} d v^{(1)} \cdots d v^{(k)},
\end{align*}
where we have set
\begin{align*}
& \bar{I}_{l_0, l_1, \cdots, l_{k-1},l_{k}}^{j_1, \cdots, j_k}   = \int_{O_N(\R)^k} (\Pi_{r=1}^k O_{l_{r-1},j_r}^{(r)} O_{l_{r},j_r}^{(r)} )   \\
& \times e^{ \frac{1}{(1+c)}\sum_{r<l} \sigma_{|x_r-x_l|}^2 \sum_{m_1,m_2=1}^N \beta^{j_r,j_l}_{m_1,m_2}\sum_{n_1,n_2=1}^N O_{n_1,m_1}^{(r)} O_{n_2,m_1}^{(r)} O_{n_1,m_2}^{(l)} O_{n_2,m_2}^{(l)}  }  dO^{(1)} \cdots dO^{(k)}  
\end{align*}
with $\beta^{j_r,j_l}_{m_1,m_2}=(-c+(1+c)1_{m_1=j_r})(-c+(1+c)1_{m_2=j_l})$. This leads to
\begin{align*}
& \bar{I}_{l_0, l_1, \cdots, l_{k-1},l_{k}}^{j_1, \cdots, j_k}  = \int_{O_N(\R)^k} (\Pi_{r=1}^k O_{l_{r-1},j_r}^{(r)} O_{l_{r},j_r}^{(r)} )   \\
& \times e^{ \frac{1}{(1+c)}\sum_{r<l} \sigma_{|x_r-x_l|}^2  (  c^2 tr (O^{(r)}{}^t O^{(r)} O^{(l)} {}^t O^{(l)}   )  +(1+c)^2 ({}^t O^{(r)} O^{(l)})_{j_r,j_l} ^2)   }  \\
& \times e^{ -c \sum_{r<l} \sigma_{|x_r-x_l|}^2  ( ({}^t O^{(r)} O^{(l)} {}^tO^{(l)} O^{(r)})_{j_r,j_r}+ ({}^t O^{(l)} O^{(r)} {}^tO^{(r)} O^{(l)})_{j_l,j_l}  )   }  dO^{(1)} \cdots dO^{(k)}  \\
& =e^{ (\frac{c^2N}{(1+c)}-2c)\sum_{r<l} \sigma_{|x_r-x_l|}^2}   \int_{O_N(\R)^k} (\Pi_{r=1}^k O_{l_{r-1},j_r}^{(r)} O_{l_{r},j_r}^{(r)} )    \\
&\times e^{(1+c) \sum_{r<l} \sigma_{|x_r-x_l|}^2   ({}^t O^{(r)} O^{(l)})_{j_r,j_l} ^2  }    dO^{(1)} \cdots dO^{(k)}  \\
\end{align*}

In conclusion, we get the following equivalent
\begin{align*}
\bar{Z} &   \underset{\epsilon \to 0}{\sim}    
\frac{\sigma_{\epsilon}^{N(N+1)k/2} e^{k\sigma_{\epsilon}^2/2} (1+c)^{(N-1)k}\sigma_{\epsilon}^{(N-1)k}}{R_N^{\epsilon}} \int_{\R^{kN}} \Pi_{r=1}^k |\Delta_{1}(v^{(r)})|  e^{- \frac {1} {2 (1+c)} \sum_{r=1}^k \sum_{i=1}^N (v_i^{(r)})^2  }   \\
& \times  e^{ - \alpha \sum_{r=1}^k  (\sum_{i=1}^N v_i^{(r)})^2 }  d v^{(1)} \cdots d v^{(k)} e^{ -c \sum_{r<l} \sigma_{|x_r-x_l|}^2} F_{i,j}(x_1, \cdots, x_k) 
 \end{align*}
where we have the following expression for $F_{i,j}(x_1, \cdots, x_k)$
\begin{align*}
 & F_{i,j}(x_1, \cdots, x_k)  =  \\
 &  \sum_{j_1, \cdots, j_k=1}^{N} \sum_{\underset{l_0=i; l_k=j}{l_1, \cdots, l_{k-1}=1}}^N \int_{O_N(\R)^k} (\Pi_{r=1}^k O_{l_{r-1},j_r}^{(r)} O_{l_{r},j_r}^{(r)} )    e^{(1+c) \sum_{r<l} \sigma_{|x_r-x_l|}^2   ({}^t O^{(r)} O^{(l)})_{j_r,j_l} ^2  }    dO^{(1)} \cdots dO^{(k)} \\
& =  \sum_{j_1, \cdots, j_k=1}^{N}  \int_{O_N(\R)^k}     O_{i,j_1}^{(1)}   ({}^tO^{(1)}O^{(2)})_{j_1,j_2} \cdots ({}^tO^{(k-1)}O^{(k)})_{j_{k-1},j_k} O_{j,j_k}^{(k)}   \\
&\times    e^{(1+c) \sum_{r<l} \sigma_{|x_r-x_l|}^2   ({}^t O^{(r)} O^{(l)})_{j_r,j_l} ^2  }    dO^{(1)} \cdots dO^{(k)}
\end{align*}

We get thus the following expression
\begin{equation}  \label{equivbarZ}
\bar{Z} /(Z c_\epsilon^k)   \underset{\epsilon \to 0}{\sim}  C_{k,N}  e^{ -c \sum_{r<l} \sigma_{|x_r-x_l|}^2} F_{i,j}(x_1, \cdots, x_k), 
 \end{equation}
where $C_{k,N}$ is a constant which depends only on $k,N$ (we can compute this constant but it is tedious and will not be necessary for the purpose of this paper).

\subsubsection{Computation of the moment of order $k$ and of the structure functions}
From relation (\ref{equivbarZ}), we get the following expression
\begin{equation}\label{orderk}
E  [  {\rm tr} M( A  )^k    ]=\underset{\epsilon \to 0}{\lim} \: E({\rm tr} M^{\epsilon}(A)^k)  = C_{k,N}  \int_{A^k}   e^{ -c \sum_{r<l} \sigma_{|x_r-x_l|}^2} \sum_{i=1}^N F_{i,i}(x_1, \cdots, x_k)   dx_1 \cdots dx_k  
\end{equation}
The main difficulty is to study the functions $F_{i,j}$. If we take the trace, we get
\begin{align}
& \sum_{i=1}^N F_{i,i}(x_1, \cdots, x_k)  \nonumber \\
& =  \sum_{j_1, \cdots, j_k=1}^{N}  \int_{O_N(\R)^k}   ({}^tO^{(1)}O^{(2)})_{j_1,j_2} \cdots ({}^tO^{(k-1)}O^{(k)})_{j_{k-1},j_k} ({}^tO^{(k)}O^{(1)})_{j_{k},j_1} \nonumber  \\  
& \times  e^{(1+c) \sum_{r<l} \sigma_{|x_r-x_l|}^2   ({}^t O^{(r)} O^{(l)})_{j_r,j_l} ^2  }    dO^{(1)} \cdots dO^{(k)}. \label{devel}
\end{align}
In particular, for $k=2$, we recover that
\begin{equation*}
\sum_{i=1}^N F_{i,i}(x_1,x_2) = N^2  \int_{O_N(\R)^2}   ({}^tO^{(1)}O^{(2)})_{1,1}^2   e^{(1+c) \sigma_{|x_2-x_1|}^2   ({}^t O^{(1)} O^{(2)})_{1,1} ^2  }    dO^{(1)} dO^{(2)}.
\end{equation*}

Here we suppose that $L=1$ and $m=0$ to simplify the presentation. Since for all $r<l$, we have $({}^t O^{(r)} O^{(l)})_{j_r,j_l} ^2 \leq 1$, it is easy to see that
\begin{align*}
e^{ -c \sum_{r<l} \sigma_{|x_r-x_l|}^2}  \sum_{i=1}^N F_{i,i}(x_1, \cdots, x_k) & \leq N^k e^{\sum_{r<l} \sigma_{|x_r-x_l|}^2   } \\
& \leq N^k \Pi_{r<l} \frac{1}{|x_r-x_l|^{\gamma^2}} \\
\end{align*}
In view of (\ref{orderk}), this entails that
\begin{equation}\label{eq:supbound}
\overline{\underset{\ell \to 0}{\lim}} \: \frac{\ln  E  [  {\rm tr} M( B(0,\ell)  )^k    ] }{\ln \frac{1}{\ell}}  \leq  - \zeta(k).
\end{equation}

We have
\begin{align*}
E  [  {\rm tr} M( B(0,\ell)  )^k    ]  & =  C_{k,N}  \int_{B(0,\ell) ^k}   e^{ -c \sum_{r<l} \sigma_{|x_r-x_l|}^2} \sum_{i=1}^N F_{i,i}(x_1, \cdots, x_k)   dx_1 \cdots dx_k   \\
& =    \ell^{dk+c \gamma^2 \frac{k(k-1)}{2}}  \int_{B(0,1) ^k}   e^{ -c \sum_{r<l} \sigma_{|u_r-u_l|}^2} \sum_{i=1}^N F_{i,i}(\ell u_1, \cdots,\ell u_k)   du_1 \cdots du_k  \\
& = \sum_{j_1, \cdots, j_k=1}^{N} \ell^{dk+c \gamma^2 \frac{k(k-1)}{2}} \int_{O_N(\R)^k}   ({}^tO^{(1)}O^{(2)})_{j_1,j_2} \cdots ({}^tO^{(k-1)}O^{(k)})_{j_{k-1},j_k} ({}^tO^{(k)}O^{(1)})_{j_{k},j_1} \\
& \times e^{(1+c) \gamma^2 \ln \frac{1}{\ell} \sum_{r<l}  ({}^t O^{(r)} O^{(l)})_{j_r,j_l} ^2  }   \left ( \int_{B(0,1) ^k}     e^{\sum_{r<l} \sigma_{|u_r-u_l|}^2  ((1+c) ({}^t O^{(r)} O^{(l)})_{j_r,j_l} ^2 -c) }     du_1 \cdots du_k \right ) \\
& dO^{(1)} \cdots dO^{(k)}
\end{align*}

In order to prove the other side of (\ref{secondequivalent}), we now study each term in the above sum.

We fix $(j_1, \cdots, j_k)$ and $\epsilon, \delta$ small such that $\epsilon<\delta$. We have
\begin{align*} 
&  \int_{O_N(\R)^k}   ({}^tO^{(1)}O^{(2)})_{j_1,j_2} \cdots ({}^tO^{(k-1)}O^{(k)})_{j_{k-1},j_k} ({}^tO^{(k)}O^{(1)})_{j_{k},j_1} \\
& \times e^{(1+c) \gamma^2 \ln \frac{1}{\ell} \sum_{r<l}  ({}^t O^{(r)} O^{(l)})_{j_r,j_l} ^2  }   \left ( \int_{B(0,1) ^k}     e^{\sum_{r<l} \sigma_{|u_r-u_l|}^2  ((1+c) ({}^t O^{(r)} O^{(l)})_{j_r,j_l} ^2 -c) }     du_1 \cdots du_k \right ) \\
& dO^{(1)} \cdots dO^{(k)}\\
& = A_{\epsilon}+A_{\epsilon, \delta}+ A_{\delta}
\end{align*}
where 
\begin{equation*}
A_{\epsilon}=\int_{ \sum_{r<l}  ({}^t O^{(r)} O^{(l)})_{j_r,j_l} ^2  \geq \frac{k(k-1)}{2}-\epsilon} \cdots , \; \; \; A_{\epsilon, \delta} = \int_{ \frac{k(k-1)}{2}-\delta  \leq  \sum_{r<l}  ({}^t O^{(r)} O^{(l)})_{j_r,j_l} ^2  \leq \frac{k(k-1)}{2}-\epsilon} \cdots 
\end{equation*}
and $A_{\delta}$ is the $ \sum_{r<l}  ({}^t O^{(r)} O^{(l)})_{j_r,j_l} ^2  \leq \frac{k(k-1)}{2}-\delta$ part of the integral. On the event $\sum_{r<l}  ({}^t O^{(r)} O^{(l)})_{j_r,j_l} ^2  \geq \frac{k(k-1)}{2}-\epsilon$, each $|{}^t O^{(r)} O^{(l)})_{j_r,j_l}|$ is greater or equal to $\sqrt{1-\epsilon}$. In particular, we have that $|({}^t O^{(r)} O^{(r+1)})_{j_r,j_{r+1}}| \geq \sqrt{1-\epsilon}$ for all $r \leq k-1$. Notice that  $({}^tO^{(k)}O^{(1)})_{j_{k},j_1}=   ({}^tO^{(1)}O^{(2)})_{j_1,j_2} \cdots ({}^tO^{(k-1)}O^{(k)})_{j_{k-1},j_k}+ O(\epsilon)$. Therefore, we can conclude that $({}^tO^{(1)}O^{(2)})_{j_1,j_2} \cdots ({}^tO^{(k-1)}O^{(k)})_{j_{k-1},j_k}$ and $ ({}^tO^{(k)}O^{(1)})_{j_k,j_1}$ have the same sign. Thus, we get
\begin{align*}
& A_{\epsilon} \\
& \geq  e^{(1+c) \gamma^2 \ln \frac{1}{\ell}( \frac{k(k-1)}{2}-\epsilon) }  |  B(0,1)^k  |  \\
&\times \int_{ \sum_{r<l}  ({}^t O^{(r)} O^{(l)})_{j_r,j_l} ^2  \geq \frac{k(k-1)}{2}-\epsilon}  ({}^tO^{(1)}O^{(2)})_{j_1,j_2} \cdots ({}^tO^{(k-1)}O^{(k)})_{j_{k-1},j_k} ({}^tO^{(k)}O^{(1)})_{j_{k},j_1}  \\
&    dO^{(1)} \cdots dO^{(k)}  \\
& \geq  e^{(1+c) \gamma^2 \ln \frac{1}{\ell}( \frac{k(k-1)}{2}-\epsilon) } ( (1-\epsilon)^k  +O(\epsilon))  |  B(0,1)^k  |  \P \left (  \sum_{r<l}  ({}^t O^{(r)} O^{(l)})_{j_r,j_l} ^2  \geq \frac{k(k-1)}{2}-\epsilon \right ) \\
\end{align*}
The only thing to check is that $\sum_{r<l}  ({}^t O^{(r)} O^{(l)})_{j_r,j_l} ^2  \geq \frac{k(k-1)}{2}-\epsilon$ has a positive probability but this can be seen easily by setting one chosen element of each $O^{(r)}$, say $O_{1,j_r}^{(r)}$, very close to one. The condition $k < \frac{2d}{\gamma^2}$ ensures that
\begin{equation*}  
c_k:=  \int_{B(0,1) ^k}     e^{\sum_{r<l} \sigma_{|u_r-u_l|}^2   }     du_1 \cdots du_k < \infty
\end{equation*}
Note that we have
\begin{equation*}
|A_{\epsilon, \delta}| \leq c_k  e^{(1+c) \gamma^2 \ln \frac{1}{\ell}( \frac{k(k-1)}{2}-\epsilon) }   \P \left (   \frac{k(k-1)}{2}-\delta  \leq  \sum_{r<l}  ({}^t O^{(r)} O^{(l)})_{j_r,j_l} ^2  \leq \frac{k(k-1)}{2}-\epsilon    \right )
\end{equation*}
Therefore, one can choose $\delta$ larger than $\epsilon$ such that $|A_{\epsilon, \delta}| \leq \frac{A_{\epsilon}}{2}$.

Finally, for these choices of $\epsilon, \delta$, we have
\begin{equation*}
  | A_{\delta} |  \leq  c_k e^{(1+c) \gamma^2 \ln \frac{1}{\ell}( \frac{k(k-1)}{2}-\delta) } 
  \end{equation*}

 We thus get the following
 \begin{equation*}
 \underset{\ell \to 0}{\underline{\lim}} \frac{\ln (   A_{\epsilon}+A_{\epsilon, \delta}+ A_{\delta}  )}{\ln \frac{1}{\ell}}   \geq (1+c) \gamma^2 ( \frac{k(k-1)}{2}-\epsilon)
 \end{equation*}
 
 Since this is valid for all $\epsilon$, we get that
\begin{equation}\label{eq:limitfinale}
\underset{\ell \to 0}{\underline{\lim}} \: \frac{\ln  E  [  {\rm tr} M( B(0,\ell)  )^k    ] }{\ln \frac{1}{\ell}}  \geq  - \zeta(k)
\end{equation}

The desired result  (\ref{secondequivalent}) is then a consequence of (\ref{eq:supbound}) and of (\ref{eq:limitfinale}).

\appendix
\section{Appendix}
\subsection{Discussion about the construction of kernels}\label{reminder}
 
 In this subsection, we discuss in further detail the construction of the kernel $K$ as summarized in remark \ref{rem:kernel}. In dimension $1$ and $2$, it is plain to see that
\begin{equation}\label{def:Keps0}
\ln_+\frac{L}{|x|}=\int_0^{+\infty}(t-|x|)_+\nu_L(dt)
\end{equation}
 where the measure $\nu_L$ is given by ($\delta_u$ stands for the Dirac mass at $u$):
 $$ \nu_L(dt)=\mathbf{1}_{[0,L]}(t)\frac{dt}{t^2}+\frac{1}{L}\delta_L(dt).$$ 
Hence, for every $\mu>0$, we have
$$\ln_+\frac{L}{|x|}=\frac{1}{\mu}\ln_+\frac{L^\mu}{|x|^\mu}=\int_0^{+\infty}(t-|x|^\mu)_+\nu_{L^\mu}(dt).$$
We are therefore led to consider  $\mu>0$ such that the function  $x\mapsto (1-|x|^\mu)_+$ is positive definite, the so-called Kuttner-Golubov problem  (see \cite{golubov}). 

For $d=1$, it is straightforward to check that $(1-|x|)_+$ is positive definite. We can thus consider a Gaussian process $X^\epsilon$ with covariance kernel given by
$$
K_\epsilon(x)=\gamma^2\int_{\epsilon}^{L}(t-|x|)_+\nu_L(dt).$$
Notice that 
\begin{equation}\label{def:Keps1}
\forall x\not=0,\quad \gamma^2\ln_+\frac{L}{|x|}=\lim_{\epsilon\to 0}K_\epsilon(x)
\end{equation}  and
\begin{equation}\label{def:Keps2}
\forall \epsilon<|x|\leq L,\quad K_\epsilon(x)=\gamma^2\int_{|x|}^{L}(t-|x|)_+\nu_L(dt)=\gamma^2\ln_+\frac{L}{|x|}  .\end{equation}

In dimension 2, we can use the same strategy since Pasenchenko \cite{cf:pas} proved that the mapping $x\mapsto (1-|x|^{1/2})_+$ is positive definite over $\R^2$.  We can thus consider a Gaussian process $X^\epsilon$ with covariance kernel given by
$$K_\epsilon(x)=2\gamma^2\int_{\epsilon^{1/2}}^{L^{1/2}}(t-|x|^{1/2})_+\nu_{L^{1/2}}(dt),$$ sharing the same properties \eqref{def:Keps1} and \eqref{def:Keps2}.
  
In dimension $3$, it is not known whether the mapping $x\mapsto \ln_+\frac{L}{|x|}$ admits an integral representation of the type explained above. Nevertheless it is positive definite so that we can use the convolution techniques developed in \cite{Rob}. In dimension $4$, it is not positive definite \cite{Rob} so that another construction is required. We explain the methods in \cite{RV1}. We set the dimension $d$ to be larger than $ d\geq 3$. Let us denote by  $S$ the sphere of  $\R^d$ and $ \sigma$ the surface measure on the sphere  such that $\sigma(S)=1$. Remind that this measure is invariant under rotations. We define the function
\begin{equation}\label{def:F}
 \forall x\in\R^d,\quad F(x)=\gamma^2\int_S\ln_+\frac{L}{|\langle x,s\rangle|}\sigma(ds).
\end{equation}
It is plain to see that $F$ is an isotropic function. Let us compute it over a neighborhood of $0$:
for $|x| \leq L $, we can write $x=|x|e_x$ with $e_x\in S$. Then we have
$$F(x)=\gamma^2\int_S\ln\frac{L}{|x||\langle e_x,s\rangle|}\sigma(ds)=\lambda^2\ln\frac{L}{|x|}+\int_S\ln\frac{1}{|\langle e_x,s\rangle|}\sigma(ds).$$
Notice that the integral  $  \int_S\ln\frac{1}{|\langle e_x,s\rangle|}\sigma(ds)$ is finite (use Lemma \ref{haargauss} below for instance) and does not depend on $x$ by invariance under rotations of the measure $\sigma$.  By using the decomposition \eqref {def:Keps0}, we can thus consider a Gaussian process $X^\epsilon$ with covariance kernel given by
$$K_\epsilon(x)= \gamma^2\int_S\int_{\epsilon}^{L}(t-|\langle x,s\rangle|)_+\nu_L(dt)\sigma(ds),$$ sharing the  properties
\begin{equation}
\forall x\not=0,\quad \lim_{\epsilon\to 0}K_\epsilon(x)=F(x)
\end{equation}  and
\begin{equation}
\forall \epsilon<|x|\leq L,\quad K_\epsilon(x)=F(x)=\lambda^2\ln\frac{L}{|x|}+C 
\end{equation}
 for some constant $C$.
 
\subsection{Auxiliary results}

We give a proof of the following standard result
\begin{lemma}\label{haargauss}
 If $(Z_i)_{1\leq i\leq N}$ are i.i.d. standard Gaussian random variables then the vector
$$V=\frac{1}{\sqrt{\sum_{i=1}^NZ_i^2}} (Z_1,\dots,Z_N)$$is distributed as the Haar measure on the sphere of $\R^N$. In particular, the density of the first entry of a random vector uniformly distributed on the sphere is given by: 
$$\frac{ \Gamma(\frac{N}{2})}{\Gamma(\frac{1}{2})\Gamma(\frac{N-1}{2})}y^{-\frac{1}{2}}
(1-y )^{\frac{N-3}{2}} \mathbf{1}_{[0,1]}(y)\,dy.
$$
\end{lemma}
 
 \noindent {\it Proof.} By using the invariance under rotations of the law of the Gaussian  vector $(Z_i)_{1\leq i\leq N}$, the law of $V$ is invariant under rotations and is supported by the sphere. By uniqueness of the Haar measure, $V$ is distributed as the Haar measure. We have to compute the density of  $\zeta_1=\frac{Z_1^2}{ \sum_{i=1}^NZ_i^2 }$. Notice that
 $$\zeta_1 =\frac{Y}{Y+Z}$$ where $Y,Z$ are independent random variables with their respective laws being chi-squared distributions of parameters $1$ and $N-1$. Therefore
 \begin{align*}
 E[f(\zeta_1 )]&=  \int_{\R_+}\int_{\R_+}f\Big(\frac{x}{x+y}\Big) \frac{1}{2^{\frac{1}{2}}\Gamma(\frac{1}{2})}x^{-\frac{1}{2}}e^{-\frac{x}{2}} \frac{1}{2^{\frac{N-1}{2}}\Gamma(\frac{N-1}{2})}y^{\frac{N-1}{2}}e^{-\frac{y}{2}}\,dx\,dy\\
 &= \frac{1}{2^{\frac{N}{2}}\Gamma(\frac{1}{2})\Gamma(\frac{N-1}{2})}\int_{0}^1f(u)\frac{1}{\sqrt{u}(1-u)^{\frac{3}{2}}}\int_{\R_+}e^{-\frac{y}{2(1-u)}}y^{\frac{N-2}{2}}\,dy\,du\\
 &=\frac{\Gamma(\frac{N}{2})}{\Gamma(\frac{1}{2})\Gamma(\frac{N-1}{2})}\int_{0}^1f(u)u^{-\frac{1}{2}}(1-u)^{\frac{N-3}{2}} \,du.
 \end{align*}
\qed
 
Next we characterize all the symmetric Gaussian random matrices
\begin{lemma}\label{isotgauss}
Let $X$ be a symmetric and isotropic centered Gaussian random matrix of size $N\times N$. Then the diagonal terms $(X_{11},\dots,X_{NN})$ have a covariance matrix of the form $\sigma^2(1+c)\mathrm{I}_N-c\sigma^2P$ for some $\sigma^2\geq 0$ and $c\in]-1,\frac{1}{N-1}]$, where $P$ is the $N\times N$ matrix whose all entries are $1$. The 
off-diagonal terms are i.i.d with variance $\sigma^2\frac{1+c}{2}$ and are independent of the diagonal terms.
\end{lemma}

 \noindent {\it Proof.} If  $X$ admits a density with respect to the Lebesgue measure $dM$ over the set of symmetric matrices (see \cite[chapter 4]{Anderson}), then the density of $M$ is given by
$$e^{-f(M)} \,dM,$$ where $f$ is a  homogeneous polynomial of degree $2$. By isotropy, $f$ must be a symmetric function of the eigenvalues of $M$. Therefore it takes on the form
$$f(M)=\alpha \mathrm{tr}(M^2)+\beta \mathrm{tr}(M)^2 $$ for some $\alpha,\beta\in\R$. In this case, the result follows easily.

If  $M$ does not admit a density with respect to the Lebesgue measure over the set of symmetric matrices, we can add an independent ``small GOE", i.e. we consider $M+\epsilon M'$ where $M'$ is a matrix of the GOE ensemble with a normalized variance independent of $M$. The matrix $M+\epsilon M'$ admits a density so that we can apply the above result. Then we pass to the limit as $\epsilon\to 0$.\qed
 
\subsection{Some integral formulae}
 Let $\alpha, c >0$. We want to compute the integral
 \begin{equation*}
\int_{\R^N}  e^{-\alpha(\sum_{i=1}^{N} \lambda_{i})^2-\frac{1}{2(1+c)}\sum_{i=1}^{N} \lambda_{i}^2}   \Pi_{i<j} |\lambda_j-\lambda_i |  d \lambda.
 \end{equation*}
 We write the integrand in the form (\ref{eq:densityvp}):
\begin{equation*}
 e^{-\alpha(\sum_{i=1}^{N} \lambda_{i})^2-\frac{1}{2\sigma_d^2(1+\bar{c})}\sum_{i=1}^{N} \lambda_{i}^2} \Pi_{i<j} |\lambda_j-\lambda_i|
\end{equation*}
where $\sigma_d^2(1+\bar{c})=(1+c)$ and $\alpha=\frac{\bar{c}}{2 \sigma_d^2(1+\bar{c})}\frac{1}{(1+\bar{c}(1-N))}$. In that case, we have $\bar{c}=\frac{2 \alpha(1+c)}{1+2 \alpha(1+c)(N-1)}$ and $1+\bar{c}(1-N)= \frac{1}{1+2 \alpha (1+c)(N-1)} $. We deduce
\begin{equation}\label{eq:integral1}
\int_{\R^N}  e^{-\alpha(\sum_{i=1}^{N} \lambda_{i})^2-\frac{1}{2(1+c)}\sum_{i=1}^{N} \lambda_{i}^2}   \Pi_{i<j} |\lambda_j-\lambda_i |  d \lambda=N! (2\pi)^{N/2} (\prod_{k=1}^N \frac{\Gamma(k/2)}{\Gamma(1/2)}) \frac{(1+c)^{N(N+1)/4}}{\sqrt{ 1+2 \alpha(1+c)N }}
\end{equation}

We also want to compute the integral
 \begin{equation*}
\int_{\R^N}  e^{-\alpha(\sum_{i=1}^{N} \lambda_{i})^2-\frac{1}{2(1+c)}\sum_{i=1}^{N} \lambda_{i}^2}   \Pi_{2 \leq i<j} |\lambda_j-\lambda_i |  d \lambda.
 \end{equation*}

We have
\begin{align*}
& \int_{\R^N}  e^{-\alpha (\sum_{i=1}^{N} \lambda_{i})^2-\frac{1}{2(1+c)}\sum_{i=1}^{N} \lambda_{i}^2}  \underset{2 \leq i<j}{\Pi} |\lambda_j-\lambda_i| d \lambda_1 \cdots d \lambda_N    \\
& = \int_{\R^{N-1}} ( \int_{\R}   e^{-2 \alpha \lambda_1 (\sum_{i=2}^{N} \lambda_{i})-(\alpha +\frac{1}{2(1+c)}) \lambda_{1}^2}   d \lambda_1)  e^{-\alpha (\sum_{i=2}^{N} \lambda_{i})^2-\frac{1}{2(1+c)}\sum_{i=2}^{N} \lambda_{i}^2}  \underset{2 \leq i<j}{\Pi} |\lambda_j-\lambda_i| d \lambda_2 \cdots d \lambda_N    \\ 
& = \sqrt{ 2 \pi} \sqrt{ \frac{1+c}{2 \alpha(1+c)+1}} \int_{\R^{N-1}}  e^{ 2 \alpha^2 \frac{1+c}{2 \alpha(1+c)+1}(\sum_{i=2}^{N} \lambda_{i})^2   -\alpha (\sum_{i=2}^{N} \lambda_{i})^2-\frac{1}{2(1+c)}\sum_{i=2}^{N} \lambda_{i}^2}  \underset{2 \leq i<j}{\Pi} |\lambda_j-\lambda_i| d\lambda_2 \cdots d \lambda_N    \\ 
& = \sqrt{ 2 \pi} \sqrt{ \frac{1+c}{2 \alpha(1+c)+1}} \int_{\R^{N-1}}  e^{ -\frac{\alpha}{2 \alpha(1+c)+1}(\sum_{i=2}^{N} \lambda_{i})^2   -\frac{1}{2(1+c)}\sum_{i=2}^{N} \lambda_{i}^2}  \underset{2 \leq i<j}{\Pi} |\lambda_j-\lambda_i| d \lambda_2 \cdots d \lambda_N    \\ 
& = \sqrt{ 2 \pi} \sqrt{ \frac{1+c}{2 \alpha(1+c)+1}} \int_{\R^{N-1}} e^{ -\frac{\bar{c}}{2 \sigma_d^2(1+\bar{c})}\frac{1}{(1+\bar{c}(2-N))}(\sum_{i=2}^{N} \lambda_{i})^2 -\frac{1}{2\sigma_d^2(1+\bar{c})}  \sum_{i=2}^{N} \lambda_{i}^2 } \underset{2 \leq i<j}{\Pi} |\lambda_j-\lambda_i| d \lambda_2 \cdots d \lambda_N    \\ 
\end{align*}
for $\sigma_d^2(1+\bar{c})=1+c$ and $\bar{c}=\frac{2 \alpha(1+c)}{2 \alpha(1+c)(N-1)+1}$ (or equivalently, $1+\bar{c}(2-N)=\frac{1+2\alpha(1+c)}{1+2\alpha(1+c)(N-1)}$ and $1+\bar{c}=\frac{1+2\alpha(1+c)N}{1+2\alpha(1+c)(N-1)}$). This leads to the following
\begin{equation*}
\int_{\R^{N-1}} e^{ -\frac{\bar{c}}{2 \sigma_d^2(1+\bar{c})}\frac{1}{(1+\bar{c}(2-N))}(\sum_{i=2}^{N} \lambda_{i})^2 -\frac{1}{2\sigma_d^2(1+\bar{c})}  \sum_{i=2}^{N} \lambda_{i}^2 } \underset{2 \leq i<j}{\Pi} |\lambda_j-\lambda_i| d \lambda_2 \cdots d \lambda_N = \bar{Z}_{N-1}
\end{equation*}
where $\bar{Z}_{N-1}=(N-1)! (2\pi)^{(N-1)/2} (\prod_{k=1}^{N-1} \frac{\Gamma(k/2)}{\Gamma(1/2)}) \sigma_{d}^{N(N-1)/2} (1+\bar{c})^{(N-2)(N+1)/4}\sqrt{1+\bar{c}(2-N)}$. In conclusion, we get:
\begin{align}
 & \int_{\R^N}  e^{-\alpha (\sum_{i=1}^{N} \lambda_{i})^2-\frac{1}{2(1+c)}\sum_{i=1}^{N} \lambda_{i}^2}  \underset{2 \leq i<j}{\Pi} |\lambda_j-\lambda_i| d \lambda_1 \cdots d \lambda_N   \nonumber \\
 & = \sqrt{1+c} (N-1)! (2\pi)^{N/2} (\prod_{k=1}^{N-1} \frac{\Gamma(k/2)}{\Gamma(1/2)})\frac{(1+c)^{N(N-1)/4}}{\sqrt{1+2 \alpha (1+c)N}}  \label{eq:integral2}   \\ \nonumber
\end{align}

\hspace{0.1 cm}

 \subsection{Heuristic derivation of the conjecture}
 Let $\ell<1$. 
We can roughly write as $\ell \to 0$ (where $\approx$ means equivalent to a random constant of order $1$)
\begin{equation*}
  M( B(0,\ell)  ) \approx \ell^d \frac{e^{ \gamma \sqrt{ \ln \frac{1}{\ell} } \Omega - \frac{\gamma^2}{2} \ln \frac{1}{\ell}}}{ \gamma^{N-1}  (\ln \frac{1}{\ell})^{(N-1)/2}    },
\end{equation*}
where $\Omega$ is a random matrix whose density is given by (\ref{eq:1pointDensity}) with $\sigma_d^2=1$, and thus we get (we forget terms of order 1)
\begin{align*} 
E  [  {\rm tr} M( B(0,\ell)  )^q    ]  & \approx  \frac{\ell^{  (d + \frac{\gamma^2}{2})q }}{ (\ln \frac{1}{\ell})^{q(N-1)/2}  }   E  [  {\rm tr} e^{ \gamma q \sqrt{ \ln \frac{1}{\ell} } \Omega}    ]   \\
& \approx   \frac{\ell^{  (d + \frac{\gamma^2}{2})q }}{ (\ln \frac{1}{\ell})^{q(N-1)/2}  } \int_{\R^N} e^{ \gamma q \sqrt{ \ln \frac{1}{\ell} } u_1}e^{-\alpha (\sum_{i=1}^{N} u_{i})^2-\frac{1}{2(1+c)}\sum_{i=1}^{N} u_{i}^2} \Pi_{i<j} |u_j-u_i| d u_1 \cdots d u_N,
\end{align*}
where $\alpha=\frac{c}{2 (1+c)}\frac{1}{(1+c(1-N))}$
Thus, if $q>0$
\begin{equation*}
E  [  {\rm tr} M( B(0,\ell)  )^q    ]   \approx \frac{\ell^{  (d + \frac{\gamma^2}{2})q }}{ (\ln \frac{1}{\ell})^{(q-1)(N-1)/2}  }
\end{equation*}
\subsection*{Acknowledgements}

The authors wish to thank Krzysztof Gaw\c{e}dzki, Alice Guionnet and Raoul Robert for fruitful discussions, and grant ANR-11-JCJC CHAMU for financial support.


\end{document}